\documentclass[12pt,letterpaper,titlepage]{amsart}
\usepackage{amsmath, amssymb, amsthm, amsfonts,amscd,xr}
\input epsf.tex
\newcommand{\N}{\mathbb{N}}

\newcommand{\R}{\mathbb{R}}
\newcommand{\C}{\mathbb{C}}

\def\beq{\begin{equation}}
\def\eeq{\end{equation}}

\def\arr{\hbox to 20pt{\rightarrowfill}}

\def\la{\langle}
\def\ra{\rangle}

\def\1{{\bf 1}}

\def\e{\varepsilon}
\def\Hom{\operatorname{Hom}}
\def\End{\operatorname{End}}
\def\tr{\operatorname{tr}}

\def\A{\mathcal A}

\def\H{\mathcal H}
\def\O{\mathcal O}

\def\cR{\mathcal R}

\def\U{\mathcal U}

\def\bfg{\mathbf g}
\def\dist{\mathbf d}
\def\mg{{\mathrm{min}}}
\def\rC{\mathrm {C}}

\newenvironment{res}
               {\begin{equation}
\begin{minipage}{0.85\textwidth}}
               { \end{minipage}\end{equation} }
\def\ber{\begin{res} }
\def\eer{\end{res}}

\numberwithin{equation}{section}
\newtheorem{thm}{Theorem}[section]

\newtheorem{lemma}[thm]{Lemma}
\newtheorem{lem}[thm]{Lemma}

\newtheorem{cor}[thm]{Corollary}
\newtheorem{prop}[thm]{Proposition}

\newtheorem{dfn}[thm]{Definition}

\theoremstyle{definition}
\newtheorem{rem}[thm]{Remark}
\newtheorem{ex}[thm]{Example}

\makeatletter
\def\section{\@startsection {section}{1}{\z@}{3.5ex plus 1ex minus
    .2ex}{2.3ex plus .2ex}{\large\bf}}
    \def\subsection{\@startsection{subsection}{2}{\z@}{3.25ex plus 1ex minus
 .2ex}{1.5ex plus .2ex}{\bf}}
\makeatother
\def\pf{{\em Proof}.\, }

\def\Ad{\operatorname{Ad}}

\def\e{\epsilon}

\def\af{\mathfrak{a}}

\def\gf{\mathfrak{g}}

\def\kf{\mathfrak{k}}

\def\nf{\mathfrak{n}}

\def\pf{\mathfrak{p}}

\def\tf{\mathfrak{t}}

\def\v{{\bf v}}
\def\w{{\bf w}}
\def\0{{\bf 0}}

\def\A{\mathcal{A}}

\def\Ad{\operatorname{Ad}}
\def\ad{\operatorname{ad}}

\def\O{\mathcal{O}}

\def\H{\mathcal{H}}
\def\f{{\bf f}}

\def\Res{\operatorname{Res}}

\def\Span{\operatorname{span}}

\def\bs{\backslash}

\makeindex
\begin{document}
\title[Analytic representations]
{Analytic representation theory  of Lie groups: General theory and
analytic globalizations of Harish--Chandra modules}
\begin{abstract} In this article a general framework for studying
analytic representations of a real Lie group $G$ is introduced.
Fundamental topological properties of the representations are analyzed.
A notion of temperedness for analytic representations is introduced,
which indicates the existence of an action of a certain natural algebra
$\A(G)$ of analytic functions of rapid decay.
For reductive groups every Harish-Chandra module $V$
is shown to admit a unique tempered analytic globalization,
which is generated by $V$ and $\A(G)$ and which embeds as the
space of analytic vectors in all Banach globalizations of $V$.
\end{abstract}
\author[Gimperlein]{Heiko Gimperlein}
\author[Kr\"otz]{Bernhard Kr\"otz}
\author[Schlichtkrull]{Henrik Schlichtkrull}

\address{Leibniz Universit\"at Hannover\\Institut f\"ur Analysis \\
Welfengarten 1\\ D-30167  Hannover
\\email: gimperlein@math.uni-hannover.de}
\address{Leibniz Universit\"at Hannover\\Institut f\"ur Analysis \\
Welfengarten 1\\ D-30167  Hannover
\\email: kroetz@math.uni-hannover.de}
\address{Department of Mathematics\\
University of Copenhagen\\Universitetsparken 5 \\
DK-2100 Copenhagen
\\email: schlicht@math.ku.dk}

\date{\today}
\thanks{}
\maketitle

\section{Introduction}

While analytic vectors are basic objects in the representation
theory of real Lie groups, a coherent framework to study general
analytic representations has been lacking so far. It is the aim of
this article to introduce categories of tempered and non--tempered
such representations and to analyze their fundamental properties.
For a representation $(\pi, E)$ of a Lie group $G$
on a locally convex space $E$ to be
\emph{analytic}, we are going to require that every vector in $E$ be
analytic and that the topology on the space of analytic vectors
coincide with the topology of $E$. No completeness assumptions on
$E$ are imposed, so that the quotient of an analytic representation
by a closed invariant subspace is again analytic.

Recall that a vector $v \in E$ is called analytic provided that the
orbit map $\gamma_v: x \mapsto \pi(x)v$ extends to a holomorphic
$E$--valued function in a neighborhood of $G$ within the
complexification $G_\C$. The space $E^\omega$ of analytic vectors
carries a natural inductive limit topology $E^\omega = \lim_{n \to
\infty} E_n$,
$$E_n=\{ v\in E\mid \gamma_v \ \hbox{extends to a holomorphic map} \
GV_n\to E\}\ ,$$ indexed by a neighborhood basis $\{V_n\}_{n\in\N}$
of the identity in $G_\C$. The induced representation $(\pi,
E^\omega)$ turns out to be continuous and indeed satisfies
$E^\omega=(E^\omega)^\omega$ in the sense of topological vector
spaces. Every analytic representation is obtained in this way. Due
to the inductive limit structure of $E^\omega$, interesting examples
tend to involve complicated and possibly incomplete topologies. For
instance, infinite dimensional Fr\'{e}chet spaces do not carry any
irreducible analytic representations of a reductive group. Still, in
spite of examples by Grothendieck and others which show how
incomplete spaces may naturally occur, important special cases
are better behaved, like for instance
the analytic vectors associated to a Banach representation, the
algebra $\mathcal{A}(G)$ below, or the analytic globalization of a
Harish--Chandra module.

Moderately growing analytic representations allow for an additional
action by an algebra of superexponentially decaying functions. To be
specific, consider a Banach representation $(\pi, E)$. Fix a
left--invariant Riemannian metric on $G$ and let $\dist$ be the
associated distance function. The continuous functions
on $G$ decaying faster than $e^{-n\dist(\cdot,\1)}$
for all $n\in \N$ form a convolution algebra $\mathcal{R}(G)$,
which is a $G$--module
under the left regular representation. If we denote the space of
analytic vectors of $\mathcal{R}(G)$ by $\mathcal{A}(G)$, the map
\begin{equation}\label{Pidef}
\Pi : \mathcal{A}(G) \to \mathrm{End}(E^\omega), \,\, \Pi(f)v =
\int_G f(x)\ \pi(x)v\ dx\ , \end{equation} gives rise to a
continuous algebra action on $E^\omega$. More general
representations will be called $\mathcal{A}(G)$\emph{--tempered}, or
of \emph{moderate growth}, provided that the integral in
(\ref{Pidef}) converges and defines a continuous action of
$\mathcal{A}(G)$.
\par Let us now specify to the case where $G$ is a real reductive group,
and let us recall that to each admissible $G$-representation $E$
of finite length one can associate the Harish-Chandra module $E_K$
of its $K$-finite vectors. Conversely, a {\it globalization}
of a given Harish-Chandra module $V$ is an admissible representation of
$G$ with $V=E_K$. The main result for this case is now as follows:

\begin{thm}\label{main}
Let $G$ be a real reductive group. Then
every Harish-Chandra module $V$ for $G$ admits a unique $\A(G)$-tempered
analytic globalization $V^\mg$. Moreover, $V^\mg$ has the property
$V^\mg=\Pi(\A(G))V.$
\end{thm}

It follows that $E^\omega\simeq V^\mg$ for every
$\A(G)$-tempered globalization $E$ of $V$ (in particular, for every
Banach globalization).
Let us mention the relationship to the announcements of Schmid and
of Kashiwara-Schmid in \cite{S} and \cite{KaS}, which (among others) assert
that every Harish-Chandra module admits a unique
{\it minimal globalization}, which is equivalent to $E^\omega$
for all Banach globalizations $E$. Our proof of Theorem
\ref{main} is independent of this theory. Instead, we rely
on the corresponding theory by Casselman and Wallach for smooth
globalizations, which is documented in \cite{W} and more
recently in \cite{BK}.

The theorem features a worthwhile corollary, namely:

\begin{cor} For an irreducible admissible Banach representation $(\pi, E)$ of a real reductive group $G$, the space of
analytic vectors $E^\omega$ is an algebraically simple $\A(G)$-module.
\end{cor}

\section{Banach representations and {\it F}-representations}
\label{sec::reps}

All topological vector spaces $E$ considered in this
paper are assumed to be Hausdorff and locally convex.
If $E$ is a topological vector space,
then we denote by $GL(E)$ the group of isomorphisms of $E$.

Let $G$ be a connected Lie group. By a {\it representation} of $G$ we
shall
understand a continuous action
$$G\times E \to E, \ \ (g, v)\mapsto g\cdot v\ ,$$
on some topological vector space $E$.  Each representation gives rise
to a group homomorphism
$$\pi: G\to E, \ \ g\mapsto \pi(g),\  \ \pi(g)v:=g\cdot v  \quad (v\in E)\, ,$$
and it is custom to denote the representation by the symbol $(\pi, E)$.

\par A representation $(\pi, E)$  is called a
{\it Banach representation} if $E$
is a Banach space.
We say that $(\pi, E)$ is an {\it F}-representation, if
$E$ is a Fr\'echet space for which there exists a defining
family of seminorms
$(p_n)_{n\in \N}$ such that for all $n\in \N$ the action
$$G\times (E, p_n)\to (E,p_n)$$
is continuous. Here $(E, p_n)$ refers to the vector space $E$ endowed with the
locally convex structure induced by $p_n$.

\begin{rem}\label{basicremark}
{\rm (a)} Every Banach representation is an {\it F}-representation.

{\rm (b)} Let $(\pi, E)$ be an {\it F}-representation. For each $n\in \N$
let us denote by
$\hat E_n$ the Banach completion of $(E, p_n)$, i.e. the completion of the
normed space $E/\{p_n=0\}$. The action of $G$  on $(E, p_n)$ factors to a
continuous action on the normed space  $E/\{p_n=0\}$ and thus induces a
Banach representation of $G$ on $\hat E_n$.

{\rm (c)} The left regular action of $G$ on the Fr\'echet space $C(G)$
defines a representation, but in general not an {\it F}-representation.
\end{rem}

Let $E^\infty$ denote the space of smooth vectors in $E$, that is,
the vectors $v\in E$ for which the orbit map $g\mapsto\pi(g)v$ is
smooth into $E$. Then $E^\infty\subset E$ is an invariant subspace,
and it is dense if $E$ is complete. The orbit map provides an
injection of $E^\infty$ into $C^\infty(G,E)$, from which $E^\infty$
inherits a topological vector space structure. Then $(\pi,E^\infty)$
is a representation. Furthermore, $E^\infty$ is a Fr\'{e}chet space
if $E$ is a Fr\'{e}chet space, and $(\pi,E^\infty)$ is an {\it
F}-representation if $(\pi,E)$ is an {\it F}-representation. By
definition, a {\it smooth representation} is a representation for
which $E^\infty=E$ as topological vector spaces.

\subsection{Growth of representations} \label{subsec::growth}

We call a function $w: G\to \R^+$ a {\it weight} if
\begin{itemize}
\item $w$ is locally bounded,
\item $w$ is sub-multiplicative, i.e. if
$$w(gh)\leq w(g) w(h)$$
for all $g, h\in G$.
\end{itemize}

To every Banach representation $(\pi, E)$ we associate the function
$$w_\pi(g):=\|\pi(g)\| \qquad (g\in G)\, ,$$
where $\|\,\cdot\,\|$ denotes the standard operator norm. It follows
from the uniform boundedness principle that $w_\pi$ is locally
bounded. Hence $w_\pi$ is a weight.

Sub-multiplicative functions can be dominated in a geometric way. For that
let us fix a left invariant Riemannian metric ${\mathbf g}$ on $G$.
Associated to ${\mathbf g}$ we obtain the Riemannian
distance function $\dist:G\times G \to \R_{\geq 0}$.
The distance function is left $G$-invariant and
hence is
recovered as $\dist(g,h)=d(g^{-1}h)$ from the function
$$d(g):= \dist(g, {\1}) \qquad (g\in G),$$
where $\1\in G$ is the neutral element. Notice that it follows from
the elementary properties of the metric that $d$ is compatible with
the group structure in the sense that
\begin{equation} \label{propertiesd}
d(g^{-1})=d(g) \quad\text{and}\quad d(gh)\leq d(g)+d(h)
\end{equation}
for all $g,h\in G$.
In particular, $g\mapsto e^{d(g)}$ is a weight.
Note also that the metric balls
$\{g\in G|\, d(g)\leq R\}$ in $G$ are compact (\cite{Ga}, p. 74).

If $w$ is an arbitrary weight on $G$, then there exist constants $c, C>0$
(depending on $w$) such that
\begin{equation}\label{exponentialestimateweight}
w(g) \leq C e ^{c d(g)} \qquad (g\in G),
\end{equation}
(\cite{Ga}, p. 75, Lemme 3). In particular, it follows that a Banach
representation has at most exponential growth
$$\|\pi(g)\|\leq C e^{c \,d(g)}.$$
By applying Remark \ref{basicremark}(b) we obtain for an
F-representation $(\pi,E)$
with defining seminorms $(p_n)_{n\in\N}$
that for each $n$ there exist constants
$c_n,C_n$ such that
\begin{equation}\label{Frepestimate}
p_n(\pi(g)v)\leq
C_n e^{c_n d(g)} p_n(v) \qquad (g\in G, v\in E).
\end{equation}

Finally, notice that it follows from
(\ref{exponentialestimateweight}) that if $d_1(g)=d_{{\mathbf
g}_1}(g,\1)$ is the function associated with a different choice of a
$G$-invariant metric, then $d_1$ is compatible with $d$, in the
sense that there exist constants $c,C>0$ such that
$$d_1(g)\leq c d(g)+C \qquad (g\in G)$$
(and vice--versa with $d,d_1$ interchanged).

\begin{rem} Suppose that $G$ is a real reductive group and $\|\cdot\|$ is a norm of $G$
(see \cite{W}, Sect 2.A.2). Then $\|\,\cdot\,\|$ is a weight and
hence there exist constants $c_1, C_1>0$ such that
$$\log\|g\|\leq c_1 d(g)+C_1 \qquad (g\in G).$$
Conversely, by following the proof of \cite{W}, Lemma 2.A.2.2,
one finds constants $c_2,C_2>0$ such that
$$d(g)\leq c_2 \log\|g\| +C_2 \qquad (g\in G).$$
\end{rem}

\section{Analytic representations}\label{sec::analyticreps}

Let us start by setting up some notations in order to discuss the
issue of analyticity in a convenient way.

\par Let us denote by $\gf$ the Lie algebra of $G$.
To simplify the exposition we will assume that $G\subset G_\C$,
where $G_\C$ is a complex group with Lie-algebra $\gf\otimes_\R
\C=:\gf_\C$.  We stress, however, that this assumption is not
necessary, since the use of $G_\C$ essentially only takes place
locally in neighborhoods $G$.

\par We extend the left invariant metric ${\mathbf g}$ to a left $G_\C$-invariant
metric on $G_\C$ and denote the associated distance function as before by $d$.
For every $n\in \N$ we set
$$V_n:=\{ g\in G_\C \mid d(g)<{\frac1n} \}\qquad \hbox{and} \qquad U_n:= V_n\cap G\, .$$
It is clear  that the $V_n$'s, resp. $U_n$'s, form a  base of the neighborhood filter
of $\1$ in $G_\C$, resp. $G$. Note that $V_n$ is symmetric, and that
$xy\in V_n$ for all $x,y\in V_{2n}$.

\subsection{The space of analytic vectors}
\par Let $(\pi, E)$ be a representation of $G$. For each $v\in E$
we denote by
$$\gamma_v: G\to E, \ \ x\mapsto \pi(x)v\ ,$$
the associated continuous orbit map.  We call
$v$ an {\it analytic} vector if $\gamma_v$ extends to a
holomorphic $E$-valued function (see Section \ref{AppHol})
on some open neighborhood of $G$ in $G_\C$.

If $v$ is analytic, then $\gamma_v$ is a real analytic map $G\to E$.
The converse statement, that real analyticity of the orbit map
implies the analyticity of $v$, holds under the assumption that $E$
is sequentially complete. Hence our definition agrees with the
standard notion of analytic vectors for Banach representations, see
for example \cite{Nelson}, \cite{Ga}, \cite{Goodman}.

\begin{rem} \label{weak anal}
If $E$ is a Banach space
or more generally a complete {\it DF}-space (see \cite{MV}, Ch.\ 25),
then it follows from \cite{Neto} Thm.\ 1,
that $v$ is an analytic vector already if the orbit map is
{\it weakly analytic}, that is,
$\lambda\circ\gamma_v: G\to \C$ is real analytic for all $\lambda\in E'$.
Here $E'$ denotes the dual space of continuous linear forms.
\end{rem}

The space of analytic vectors is denoted by $E^\omega$.
A theorem of Nelson (\cite{Nelson} p.\ 599)
asserts that $E^\omega$ is dense in $E$ if $E$ is a Banach space.
More precisely, Nelson's theorem asserts the following. Let
$h_t\in C^\infty(G)$ denote the {\it heat kernel} on $G$,
where $t>0$, then $\Pi(h_t)v\in E^\omega$ and
$\Pi(h_t)v\to v$ for $t\to 0$, for all $v\in E$. In
fact, the proof of Nelson's theorem is valid more generally if $E$
is sequentially complete and with suitably restricted growth 
of $\pi$. In particular, this is the case for {\it F}-representations,
see (\ref{Frepestimate}). The density is false in general, as easy
examples such as the left regular representation of $\R$ on
$C_c(\R)$ show.

We wish to emphasize that $E^\omega$ is a $G$-invariant vector subspace of $E$.
This follows immediately from the identity
$\gamma_{\pi(g)v}(x)=\gamma_v(xg)$. We also note that $E^\omega$
is a $\gf$-invariant subset of the space $E^\infty$ of smooth vectors.

It is convenient to introduce the following notation.
For every $n\in \N$ we define the subspace of $E^\omega$,
$$E_n=\{ v\in E\mid \gamma_v \
\hbox{extends to a holomorphic map} \ GV_n\to E\}.$$ Since $G$ is
totally real in $G_\C$ and $GV_n$ is connected, the holomorphic
extension of $\gamma_v$ is unique if it exists. Let us denote the
extension by $\gamma_{v,n}\in\O(GV_n,E)$. For each $z\in GV_n$ the
operator
$$\pi_n(z): E_n\to E, \quad \pi_n(z)v:= \gamma_{v,n}(z),$$
is linear. In particular, uniqueness implies
$$\pi_{n}(gz)=\pi(g)\pi_{n}(z)$$
for all $g\in G$, $z\in GV_n$. It is easily seen that if $m<n$, then
$E_m\subset E_n$ and $\pi_m(z)v=\pi_n(z)v$ for $z\in GV_n$, $v\in
E_m$. We shall omit the subscript $n$ from the operator $\pi_n(z)$
if no confusion is possible.

A closely related space is
$$\tilde E_n=\{ v\in E\mid
\gamma_v|_{U_n} \hbox{extends holomorphically to} \ V_n\}.$$

\begin{lemma}\label{union En}
The space of analytic vectors is given by the increasing unions
$$E^\omega=\bigcup_{n\in \N} E_n=\bigcup_{n\in \N} \tilde E_n.$$
Furthermore,
\begin{equation}\label{Enchain}
E_n\subset \tilde E_n\subset E_{4n}
\end{equation}
for all $n\in\N$.
\end{lemma}

\begin{proof}
The inclusions
$$\bigcup_{n\in \N} E_n\subset E^\omega\subset \bigcup_{n\in \N} \tilde E_n,$$
as well as the first inclusion in (\ref{Enchain}),
are clear. Hence it suffices to prove the second inclusion in
(\ref{Enchain}).
Let $v\in V_{n}$ and let us denote the extension of $\gamma_v$
by $f: V_n\to E$. For $g\in G$ and $z\in V_{4n}$ we define
$$F(gz):=\pi(g)f(z)\in E.$$
We need to show that the expression is well-defined.
Assume $gz=g'z'$ with
$g,g'\in G$ and $z,z'\in V_{4n}$. Then $g^{-1}g'=zz'^{-1}\in V_{2n}$,
and hence $g^{-1}g' x\in V_n$ for all $x\in V_{2n}$.
Since $\pi(g)\pi(g^{-1}g' x)v=\pi(g')\pi(x)v$ for $x\in G$,
analytic continuation from $U_{2n}$ implies
$\pi(g)f(g^{-1}g' x)=\pi(g')f(x)$
for $x\in V_{2n}$. In particular, with $x=z'$ we obtain
$\pi(g)f(z)=\pi(g')f(z')$, showing that
$F$ is well defined. As $F$
is clearly holomorphic, we conclude that
$v\in E_{4n}$.
\end{proof}

Next we want to topologize $E^\omega$. For that we notice that
the holomorphic extensions provide injections of
$E_n$ and $\tilde E_n$ into $\O(GV_n,E)$ and $\O(V_n,E)$,
respectively.  We topologize $E_n$ and $\tilde E_n$
by means of these maps and the standard compact open topologies.
It is easily seen that
the inclusion maps $E_n\to E_{n+1}\to E$ and $\tilde E_n\to \tilde E_{n+1}\to E$
are all continuous. Furthermore:

\begin{lemma}\label{topology En}
The inclusion maps in (\ref{Enchain}) are continuous
for all $n\in\N$.
\end{lemma}

\begin{proof}
Identifying $E_n$ and $\tilde E_n$ with the corresponding
spaces of holomorphic functions, we obtain the following
neighborhood bases of $0$. In $E_n$, the members are all sets
$$W_{K,Z}:=\{ f\in E_n\mid f(K)\subset Z\},$$
where $K\subset GV_n$ is compact and $Z\subset E$
a zero neighborhood. Similarly in $\tilde E_n$, members are
$$\tilde W_{K,Z}:=\{ f\in \tilde E_n\mid f(K)\subset Z\},$$
where $K\subset V_n$ is compact and $Z\subset E$
a zero neighborhood. The continuity of the first inclusion is
then obvious.

With the mentioned identifications,
the second inclusion is given by the map $f\to F$ described in
the previous proof. Let a neighborhood $W=W_{K,O}\subset E_{4n}$
be given. Let $K'\subset V_{4n}$ be an arbitrary compact neighborhood
of $0$. By compactness of $K\subset GV_ {4n}$
we obtain a finite union $K\subset \cup g_i K'\subset GV_{4n}$.
Let $O'=\cap \pi(g_i)^{-1}(O)$, then $\tilde W=\tilde W_{K',O'}$
is an open neighborhood of $0$ in $\tilde E_n$, and
$f\in\tilde W \Rightarrow F\in W.$
\end{proof}

We endow $E^\omega$ with the
inductive limit topology of the ascending unions in Lemma \ref{union En}.
The Hausdorff property follows, since $E$ is assumed to be Hausdorff.
It follows from Lemma \ref{topology En}, that the two unions
give rise to the same topology. In symbols:
\begin{equation}\label{inductivelimit}
E^\omega=\lim_{n\to \infty} E_n=\lim_{n\to \infty} \tilde E_n\subset E\,
\end{equation}
with continuous inclusion into $E$. Since the restriction
$\O(GV_n,E)\to C^\infty(G,E)$ is continuous for all $n\in\N$,
we have $E^\omega\subset E^\infty$ with continuous inclusion.

Observe that an intertwining operator
$T\colon E\to F$ between two representations $(\pi,E), (\rho,F)$ carries
$E^\omega$ continuously into $F^\omega$. In fact,
if $v\in E_n$ with the holomorphically extended orbit map
$z\mapsto \pi(z)v$, then $Tv\in F_n$ since $z\mapsto
T\pi(z)v$ is a holomorphic extension of the orbit map
$g\mapsto \rho(g)Tv=T\pi(g)v$. It follows that $T$
maps $E_n$ continuously into $F_n$ for each $n$.

Notice that if
we define a continuous action of $G$ on $\O(GV_n, E)$ by
$$(g\cdot f)(z):= \pi(g) f(g^{-1}z) \qquad (g\in G, z\in GV_n),$$
then the image of $v\mapsto\pi_n(\cdot)v$ is the subspace $\O(GV_n,
E)^G$ of $G$-invariant functions, with inverse map given by
evaluation at $\1$. Thus $E_n$ is identified with a {\it closed}
subspace of $\O(GV_n,E)$. In particular, it follows (see \cite{J},
p.~365) that $E_n$ is complete/Fr\'echet if $E$ has this property.

Let us briefly recall the structure of the open neighborhoods of zero
in the limit $E^\omega$.
If $A$ is a subset of  some vector space, then we write $\Gamma(A)$ for
the convex hull
of $A$. Now given for each $n$ an open
$0$-neighborhood $W_n$ in $E_n$ (or $\tilde E_n$), the set
\begin{equation}\label{inductive topology}
W:=\Gamma (\bigcup_{n\in \N}  W_n)
\end{equation}
is an open convex neighborhood of $0$ in $E^\omega$.
The set of neighborhoods $W$ thus obtained
form a filter base of the $0$-neighborhoods in $E^\omega$.

\begin{prop}\label{p=1} Let $(\pi, E)$ be a representation of a Lie group on a
topological vector space $E$. Then the following assertions hold:
\begin{enumerate}
\item The action $G\times E^\omega\to E^\omega$ is continuous,
hence defines a representation $(\pi,E^\omega)$ of $G$.
\item  Each $v\in E^\omega$  is an analytic vector for $(\pi, E^\omega)$
and
$$(E^\omega)^\omega=E^\omega $$
as topological vector spaces.
\end{enumerate}
\end{prop}

\begin{proof}
In (i) it suffices to prove continuity at $(\1,v)$ for each $v\in
E^\omega$. We first prove the separate continuity of $g\mapsto
\pi(g)v \in E^\omega$. Let $v\in E_n$, and consider the $E$-valued
holomorphic extension of $g\mapsto\pi(g)v$. Since multiplication in
$G_\C$ is holomorphic and $V_{2n}\cdot V_{2n}\subset V_n$, it
follows that for each $z_1\in V_{2n}$ the element $\pi_{2n}(z_1)v$
belongs to $E_{2n}$, with the holomorphic extension
\begin{equation}\label{piproduct}
z_2\mapsto \pi_{2n}(z_2)\pi_{2n}(z_1)v:=\pi_n(z_2z_1)v \qquad
(z_1,z_2\in GV_{2n},\ v\in E_n)
\end{equation}
of the orbit map. In particular,
(\ref{piproduct}) holds for $z_1=g\in U_{2n}$.
The element $\pi(z_2g)v\in E$ depends
continuously on $g$, locally uniformly with respect to
$z_2$.
It follows that $g\mapsto \pi(g)v$ is
continuous $U_{2n}\to E_{2n}$, hence into $E^\omega$.

In order to conclude the full continuity of
(i) it now suffices to establish the following:

(*) {\it For all compact subsets $B\subset G$ the operators
$\{\pi(g)\mid g\in B\}$ form an equicontinuous
subset of $\End(E^\omega)$.}

Before proving this, we note that for every compact subset $B\subset
G$ and every $m\in\N$ there exists $n>m$ such that
$$ b^{-1}V_nb\subset V_m \quad (b\in B).$$
This follows from the continuity of the adjoint action. Then $zb\in
GV_m$ for all $z\in GV_n$, and hence $\pi(b)v\in E_n$ for all $b\in
B$, $v\in E_m$ with
\begin{equation}\label{B action}
\pi_n(z)\pi(b)v=\pi_m(zb)v.
\end{equation}

In order to prove (*) we fix a compact set $B\subset G$.
Given $m\in N$ we choose $n>m$ as above. We are going to
prove equicontinuity $B\times E_m\to E_n$. An open
neighborhood of $0$ in $E_n$ can be assumed of the form
$$(K, Z):=\{ f\in E_n\mid f(K)\subset Z\},$$
where $K\subset GV_n$ is compact and $Z\subset E$
a zero neighborhood. Then with
$K'= \cup_{b\in B} b^{-1}Kb$ and $Z'=\cap_{b\in B} \pi(b)^{-1}(Z)$
we obtain
$$f(K') \subset Z' \Rightarrow  \pi(b)f(b^{-1}Kb)\subset Z$$
for all $b\in B$ and all functions $f: GV_m\to E$.
If in addition $f$ is $G$-invariant, then
the conclusion is $f(Kb)\subset Z$, and we have
shown that the right translation by $b$ maps the zero neighborhood
$(K',Z')$ in $E_m$ into the zero neighborhood $(K,Z)$ in $E_n$.

The equicontinuity $B\times E^\omega\to E^\omega$ is an easy
consequence given the description (\ref{inductive topology}) of the
neighborhoods in the inductive limit. This completes the proof of
(i).

For the proof of (ii), let  $v\in E_n$. In the first part of the proof
we saw that $\pi(z_1)v\in E_{2n}$ for
each $z_1\in V_{2n}$, with the holomorphically extended orbit map
given by (\ref{piproduct}).
It then follows from Lemma \ref{lemma: productiso}, applied
to $V_{2n}\times GV_{2n}$ and the
map $(z_1,z_2)\mapsto \pi(z_2z_1)v$,
that $z_1\mapsto \pi(\cdot)\pi(z_1)v$ is holomorphic
$V_{2n}\to \O(GV_{2n},E)$. Hence
$z_1\mapsto \pi(z_1)v$
is holomorphic into $E_{2n}$, hence also into $E^\omega$.
Thus $g\mapsto \pi(g)v$ extends to a holomorphic $E^\omega$-valued
map on $V_{2n}$, and hence $v\in (E^\omega)^\omega$ by the
second description in (\ref{inductivelimit}).

For the topological statement in (ii), we need to show that
the identity map is continuous $E^\omega\to (E^\omega)^\omega$.
We just saw that the identity map takes
$$E_n \to \widetilde{(E^\omega)}_{2n},$$
hence it suffices to show continuity of this map
for each $n$. The proof given above reduces to the statement
that the map mentioned below (\ref{productiso}) is continuous.
\end{proof}

\begin{cor}\label{analytic-smooth}
 $(E^\infty)^\omega=(E^\omega)^\infty=E^\omega$ as topological vector spaces.
\end{cor}

\begin{proof}
The continuous inclusions
$E^\omega\subset E^\infty\subset E$ induce continuous inclusions
$E^\omega=(E^\omega)^\omega\subset (E^\infty)^\omega\subset E^\omega$.
With $E$ replaced by $E^\omega$, the same inclusions also imply
$(E^\omega)^\omega\subset (E^\omega)^\infty\subset E^\omega$.
\end{proof}

We are interested in the functorial properties of the
construction.

\begin{lemma}\label{subrep}
Let $(\pi,E)$ be a representation, and let $F\subset E$ be a closed
invariant subspace. Then
\begin{enumerate}
\item $F^\omega=E^\omega\cap F$ as a topological space,
\item $E^\omega/F^\omega \subset (E/F)^\omega$ continuously.
\end{enumerate}
\end{lemma}

\begin{proof} (i) Obviously $F_n\subset E_n$ for all $n$. Conversely,
if $v\in E_n\cap F$ with holomorphically extended
orbit map $z\mapsto\pi(z)v\in E$, then $\pi(g)v\in F$ for all $g\in G$
implies $\pi(z)v\in F$ for all $z\in GV_n$. Hence $v\in F_n$.
The topological statement follows easily.

(ii) The quotient map induces a continuous map $E^\omega\to (E/F)^\omega$,
which in view of (i) induces the mentioned continuous inclusion.
\end{proof}

Notice also that if $E_1,E_2$ are representations,
then the product representations satisfy
$E_1^\omega\times E_2^\omega\simeq (E_1\times E_2)^\omega$.

\subsection{Completeness}

In general completeness of $E$ does not ensure
that $E^\omega$  is complete.
For Banach representations this is the case as the
following result shows.

\begin{prop}\label{p=c}  Let $(\pi, E)$ be a representation of $G$
on a complete {\it DF}-space. Then $E^\omega$ is
complete.
\end{prop}

\begin{proof}
Let $(v_i)$ be a Cauchy net in $E^\omega$. It is Cauchy in $E$, hence
converges to some element $v\in E$. Moreover, the net of orbit maps
$(\gamma_{v_i})$ converges pointwise on $G$ to $\gamma_v$.
We need to show that $\gamma_v$ is real analytic, and using our
assumptions on $E$ it suffices to prove weak analyticity,
see Remark \ref{weak anal}.

Let $K\subset G$ be any compact set.
We consider the space $A(K)$ of real analytic functions on $K$.
These are germs of holomorphic functions defined on open neighborhoods
$V$ of $K$ in $G_\C$, and $A(K)$ is equipped with the inductive topology.
Since each $\O(V)$ has the Montel property, the limit is compact,
so that $A(K)$ inherits completeness from $\O(V)$.

For every $\lambda\in E'$ we consider the mapping
$$E^\omega \to A(K), \ \ E_n\ni v\mapsto \hbox{germ of } \lambda\circ \gamma_v.$$
It is clear that this is a continuous map.
It follows that $\lambda\circ\gamma_{v_i}|_K$ converges in $A(K)$, so that
$\lambda\circ\gamma_v$ is real analytic on $K$.
\end{proof}

\begin{rem}
Combining the proof above with \cite{BD}, Theorem 3,
leads to a more general result for representations on
Fr\'{e}chet spaces. In this case, $E^\omega$ is complete whenever
there is a fundamental system of seminorms $\{p_n\}_{n\in \N}$ for
the topology of $E$ such that
$$\exists n\ \forall m \geq n\ \exists k\geq m\ \exists C>0\ \forall
v \in E :\ p_m(v)^2\leq Cp_k(v)p_n(v).$$
\end{rem}

\begin{rem}
An example by Grothendieck, \cite{GreineII}, p.\ 95, may be adapted
to give an example of an incomplete space of analytic vectors.
Consider the regular representation of $G=S^1$ on the (complete)
space $E=C(S^1,\C^\N)$, where $\C^\N$ is endowed with the product
topology. The analytic vectors for this action are sequences of
functions, which extend holomorphically to a common annulus $\{z \in
\C\ |\ 1-\varepsilon<|z|<1+\varepsilon\}$ for some $\varepsilon>0$.
Being a  dense subspace of $(C(S^1)^\omega)^\N$, $E^\omega$ fails to
be complete as well as sequentially complete.
\end{rem}

\subsection{Definition of analytic representation}

Motivated by Proposition \ref{p=1} we shall give the following definition.

\begin{dfn} A representation $(\pi, E)$ is
called {\it analytic} if
$E=E^\omega$ holds as topological
vector spaces. \end{dfn}

Given a representation $(\pi, E)$, Proposition \ref{p=1} implies that
$(\pi, E^\omega)$ is an analytic representation.

\begin{lemma}
Let $(\pi,E)$ be an analytic representation, and let $F\subset E$ be
a closed invariant subspace. Then $\pi$ induces analytic
representations on both $F$ and $E/F$.
\end{lemma}

\begin{proof} This follows from Lemma \ref{subrep}. From (i) in
that lemma we infer immediately that $F^\omega=F$, and from (ii) we
then conclude that $E/F=E^\omega/F^\omega \to (E/F)^\omega$ is
continuous. The opposite inclusion is trivially valid and
continuous.
\end{proof}

\begin{ex} We consider the Fr\'{e}chet space $E:=\O(G_\C)$
with the right regular action of $G$,
$$\pi(g)f(z)=f(zg) \qquad (g\in G, z\in G_\C , f\in \O(G_\C))\, .$$
It is easy to see that $(\pi, E)$ defines a representation. Given
$v\in E$, it follows from (\ref{productiso}) that the orbit map
$\gamma_v: G\to E$ extends to a holomorphic mapping from $G_\C$ to
$E$. The same equation implies easily that $E=E^\omega$ as
topological spaces. Thus $(\pi, E)$ is analytic.
\end{ex}

\subsection{Irreducible analytic representations}

It is a natural question on which type of topological vector spaces
$E$ one can model irreducible analytic representations. The next
result shows that this class is rather restrictive.

\begin{thm} \label{thm::frechetrep} Let  $(\pi, E)$ be an irreducible
representation of a reductive group on a Fr\'{e}chet space $E$. If
$E=E^\omega$ as vector spaces, then $E$ is finite dimensional.
\end{thm}

\begin{proof} By passing to a covering group if necessary,
we may assume that $G_\C$ is simply connected.
By assumption $E^\omega=\lim E_n$ identifies with $E$ as vector spaces.
The Grothendieck factorization theorem
implies that $E=E_n$ for some $n$
(see \cite{G}, Ch. 4, Sect. 5, Thm. 1).
Hence the operator $\pi(x):=\pi_n(x)$ is defined on $E$,
for all $x\in V_n$. We shall holomorphically
extend to all $x\in G_\C$.

Let $v\in E$. By the monodromy theorem it suffices to extend
$\pi(x)v$ along all simple smooth curves starting at $\1$. Let
$\gamma: [0,1]\to G_\C$ be such a curve with $\gamma(0)=\1$. We
select finitely many open sets $U_1,\dots,U_k\subset G_\C$ which
cover the curve $\gamma([0,1])$ and points
$$x_i=\gamma(t_i), \qquad
0=t_1<\dots<t_k<1,$$ such that $\1=x_1\in U_1$ and $x_i\in U_i\cap
U_{i-1}$ for $i>1$. By choosing the sets $U_i$ sufficiently small
(and sufficiently many) we may assume that $U_i\subset V_{2n}x_i$
for each $i$ and also that the only non--empty overlaps are among
neighboring sets $U_i$ and $U_{i-1}$ (to attain these properties it
may be useful from the outset to select the sets inside a tubular
neighborhood around the curve).

In particular, $\pi(x)v$ is already defined for $x\in U_1\subset V_{2n}$.
On $U_2,\dots, U_k$ we recursively define
$$\pi(x)v=\pi(z)\pi(x_i)v, \qquad x=zx_i\in U_i\subset V_{2n}x_i,$$
where $\pi(x_i)v$ is defined in the preceding step.
Clearly this depends holomorphically on $x$.
However, in order to obtain a proper extension of $x\mapsto \pi(x)v$,
we need to verify that $\pi(x)v$
is well defined on overlaps between the $U_i$.
What we need to show is that
$$\pi(z)\pi(x_i)v=\pi(zx_i)v, \qquad zx_i\in U_i\cap U_{i-1}.$$
Let $x_i=yx_{i-1}$ where $y\in V_{2n}$.
By the recursive definition we have
$\pi(x_i)v=\pi(y)\pi(x_{i-1})v$ and
$\pi(zx_i)v=\pi(zy)\pi(x_{i-1})v$.
Then the desired identity follows
since $\pi(z)\pi(y)=\pi(zy)$ by (\ref{piproduct}).

Thus the representation extends to an irreducible holomorphic
representation of $G_\C$ (also denoted by $\pi$) .
If $U<G_C$ is a compact real form, then
the Peter-Weyl theorem implies that $\pi|_U$ is irreducible
and finite dimensional.
\end{proof}

\begin{rem}
Non--reductive groups, on the other hand, may have irreducible
analytic actions on a Fr\'{e}chet space. As an example, consider the
Schr\"{o}dinger representation of the Heisenberg group
$\mathbb{H}^n$ on the Fr\'{e}chet space
$$E=\{f \in \mathcal{O}(\C^n) \mid \forall N,M \in \N:  \sup_{x \in \R^n}\sup_{y \in (-N,N)^n} |f(x+iy)| \ e^{M|x|} < \infty\}.$$ It is irreducible as a restriction of the Schr\"{o}dinger
representation on $L^2(\R^n)$, and one readily verifies that $E =
E^\omega$.
\end{rem}

\section{The algebra of analytic superdecaying
functions}\label{sec::AG}

We define a convolution algebra of analytic functions with fast
decay. The purpose is to obtain an algebra which acts on
representations of restricted growth, such as {\it F}-representations.

\subsection{Superdecaying functions}

Let us denote by $dg$ the Riemannian measure on $G$ associated to the metric $\mathbf g$ and note that
$dg$ is a left Haar measure. It is of some relevance below that there is a constant $c>0$ such that
\begin{equation} \label{hb} \int_G  e^{-c d(g)} \ dg <\infty\end{equation}
(see \cite{Ga}, p. 75, Lemme 2).
\par We define the space of superdecaying continuous function on $G$ by
$$\cR(G):=\{ f\in C(G)\ \mid\ \forall N\in \N:\ \sup_{g\in G}  |f(g)| e^{N d(g)}<\infty\}\, $$
and equip it with the corresponding family of seminorms. Note that
$\cR(G)$ is independent of the choice of the left-invariant metric,
and that it has the following properties:

\goodbreak
\begin{prop}  \
\begin{enumerate}
\item $\cR(G)$ is a Fr\'echet space and the natural action of $G\times G$
by left-right displacements defines an {\it F}-representation,
\item $\cR(G)$ becomes a Fr\'echet algebra under convolution:
$$f*h (x)=\int_G f(y) h(y^{-1}x) \ dy $$
for $f, h\in \cR(G)$ and $x\in G$,
\item Every F-representation $(\pi, E)$ of $G$ gives rise to
a continuous algebra representation of $\cR(G)$,
$$\cR(G)\times E \to E, \ \ (f, v)\mapsto\Pi(f)v,$$
where
$$\qquad\Pi(f)v:=\int_G f(g) \pi(g) v \ dg \qquad (f\in \cR(G), v\in E) \, $$
as an $E$-valued integral.
\end{enumerate}
\end{prop}

\begin{proof} Easy. Use (\ref{propertiesd}), (\ref{Frepestimate})
and (\ref{hb}).
\end{proof}

\subsection{Analytic superdecaying functions}

We shall start with a discussion of the analytic vectors in
$\cR(G)$. Henceforth we shall view $\cR(G)$ as a $G$-module for the
left regular representation of $G$. We set $\A(G):=\cR(G)^\omega$
and equip $\A(G)$ with the corresponding vector topology. With the
notation from the preceding section we put $\A_n(G):=\cR(G)_n$ for
each $n\in \N$. Notice that $\A_n(G)$ is a Fr\'echet space for each
$n$, since $\cR(G)$ is Fr\'echet. Hence $\A(G)$ is an {\it LF}-space
(inductive limit of Fr\'echet spaces). In the Appendix we show that
$\A(G)$ is complete and reflexive.

\begin{prop} \
\begin{enumerate}
\item $\A(G)$ carries representations of $G$ by left and right
action,
\item $\A(G)$ is a subalgebra of $\cR(G)$ and convolution is
continuous $$\A(G)\times\A(G)\to\A(G).$$
\end{enumerate}
\end{prop}

\begin{proof} (i) The statement about the left action is immediate
from Proposition \ref{p=1} (i).
It is clear that $\A(G)$ is right invariant,
since every right displacement $f\mapsto R_gf$ is an intertwining operator
for the left regular representation. The continuity of the right
action follows from Lemma \ref{realization of A_n} below,
see Remark \ref{right action}.

(ii) This follows from Proposition
\ref{A action} (to be proved below) by taking $E=\cR(G)$.
\end{proof}

The next lemma gives us a concrete realization of $\A_n(G)$.

\begin{lem}\label{realization of A_n}
For all $n\in \N$, restriction to $G$ provides a topological
isomorphism of
$$\Bigg\{ f \in \O(V_n G)\Biggm|
\begin{aligned}
& \forall N>0,  \forall \Omega \subset V_n  \text{ compact}: \\
&
\sup_{g\in G, z\in \Omega}
|f(zg)| e^{N d(g)} <\infty
\end{aligned}
\Bigg\}\, $$
onto $\A_n(G)$.
Here the space above is
topologized by the seminorms mentioned in its definition.
\end{lem}

\begin{proof} Let $f\in \A_n(G)$.
Then $\gamma_f: G\to \cR(G), \ g\mapsto f(g^{-1}\cdot)$
extends to a holomorphic map
$\gamma_{f,n}: GV_n\to \cR(G)$.
As point evaluations $\cR(G)\to \C$
are continuous, it follows that
$F(z):= \gamma_{f,n}(z^{-1})(\1)$ defines a holomorphic
extension of $f$ to $V_n G$. Moreover,
$F(zg)=\gamma_{f,n}(z^{-1})(g)$ for $z\in V_n, g\in G$.
Let $N>0$ and a compact set
$\Omega\subset V_n$ be given, then
$$\sup_{g\in G, z\in \Omega}
|F(zg)| e^{N d(g)}= \sup_{z\in\Omega} p_N(\gamma_{f,n}(z^{-1}))<\infty$$
where $p_N(h)=\sup_{g\in G} |h(g)|e^{Nd(g)}$
is a defining seminorm of $\cR(G)$.
Hence $F$ belongs to the space above.
Moreover, we see that $f\mapsto F$
is an isomorphism onto its image.

Conversely, let $F$ belong to the space above and put $f:=F|_G$.
Then it is clear that $f\in\cR(G)$ (take $\Omega=\{\1\}$). We need
to show that $f\in \A_n(G)$, i.~e.~that $\gamma_f: G\to \cR(G)$
extends to a holomorphic map $GV_n \to\cR(G)$. The extension is
$z\mapsto F(z^{-1}\cdot)$, and we need to show that it is
holomorphic.

We first show that $z\mapsto F(z^{-1}\cdot)$ is continuous into
$\cR(G)$. To see this, let $z_0\in GV_n$ and $\e, N>0$ be given. We
wish to find a neighborhood $D$ of $z_0$ such that
\begin{equation}\label{estimate pN}
p_N(F(z^{-1}\cdot)-F(z_0^{-1}\cdot))<\e
\end{equation} for all $z\in D$.

Let us fix a compact neighborhood $D_0$ of $z_0$ in $GV_n$.
As
$$\sup_{g\in G, z\in D_0} |F(z^{-1}g)| e^{m d(g)} <\infty$$
for all $m>N$ we find a compact subset $K\subset G$ such that
$$  \sup_{g\in G\setminus K, z\in D_0} |F(z^{-1}g)| e^{N d(g)} <\e/2.$$
Shrinking $D_0$ to some possibly smaller neighborhood $D$ we may request that
$$  \sup_{g\in K, z\in D} |F(z^{-1}g)- F(z_0^{-1}g)| e^{N d(g)} <\e\, .$$
The required estimate \ref{estimate pN} follows.

As continuity has been verified, holomorphicity follows provided
$z\mapsto \lambda(F(z^{-1}\cdot))$ is holomorphic for $\lambda$
ranging in a subset whose linear span is weakly dense in $\cR(G)'$
(see \cite{GreineI}, p. 39, Remarque 1). A convenient such subset is
$\{\delta_g\mid g\in G\}$, and the proof is complete.
\end{proof}

\begin{rem}\label{right action} Let
$q(f):=\sup_{g\in G, z\in \Omega}
|f(zg)| e^{N d(g)}$  be a seminorm on $\A_n(G)$ as above.
Then
(\ref{propertiesd}) implies
$$q(R_xf) \leq e^{N d(x)}q(f)\quad (f\in\A_n(G))\ ,$$
for $x\in G$,
so $\A_n(G)$ is an {\it F}-representation for the right action.
\end{rem}

\subsection{Analytic vectors of F-representations}

Let $(\pi, E)$ be an {\it F}-representation of $G$, and let
$v\in E$. The map $f\mapsto \Pi(f)v$ is intertwining from
$\cR(G)$ (with left action) to $E$. Hence
$\Pi(f)v\in E_n$ for $f\in\A_n(G)$ and
$\Pi(f)v\in E^\omega$ for $f\in\A(G)$. With the preceding
characterization of $\A_n(G)$ we have
\begin{equation}\label{orbitmapPi(f)v}
\pi(z)\Pi(f)v = \int_G f(z^{-1}g)\pi(g)v\,dg
\end{equation}
for $f\in\A_n(G)$, $z\in GV_n$.

\begin{rem} In particular
$$\Pi(\A(G))E^\omega\subset E^\omega$$
for F-representations. In fact one can
show  (see \cite{GKL}) that
$$\Pi(\A(G))E^\omega=E^\omega.$$
\end{rem}

It is easily seen that the action of $\A(G)$ on $E^\omega$
is an algebra action. We shall now see that it is continuous.

\begin{prop}\label{A action}
Let $(\pi, E)$ be an F-representation. The
bilinear map $(f, v)\mapsto \Pi(f)v$ is continuous
$$\A_n(G)\times E\to E_n,$$
for every $n\in \N$. Likewise, it is continuous
$$\A(G)\times E\to E^\omega.$$
\end{prop}

Notice that
since $E^\omega$ injects continuously in $E$,
the last statement implies continuity of both
$$\A(G)\times E\to E\quad\text{and}\quad
\A(G)\times E^\omega\to E^\omega.$$

\begin{proof}
Let $n\in \N$ be fixed and let $W\subset E_n$ be an open neighborhood of $0$.
We may assume
$$W=W_{K, p}:=\{ v\in E_n\mid p(\pi(K)v)<1\}$$
with $K\subset GV_n$ compact and
$p$ a continuous seminorm on $E$ such that
$$p(\pi(g)v)\leq C e^{c d(g)} p(v) \qquad (g\in G, v\in E)$$
for some constants $c,C$
(see \ref{Frepestimate}).

Choose $N>0$ so that (cf \ref{hb})
$$C_1:=\int_G e^{(c-N)d(g)}\,dg<\infty,$$
and let
$$O:=\{ f\in\O(V_nG)\mid \sup_{z\in K, g\in G} |f(z^{-1}g)| e^{Nd(g)} <\epsilon\}
\subset \A_n(G)$$
(with $\epsilon$ to be specified below).
According to Lemma \ref{realization of A_n}, $O$ is open.

For $f\in O$ and $z\in K$ we obtain by (\ref{orbitmapPi(f)v})
\begin{equation*}
p(\pi(z)\Pi(f)v)\leq
 \int_G |f(z^{-1} g)|\, p(\pi(g)v)) \ dg \leq \epsilon CC_1 p(v) .
\end{equation*}
With $\epsilon<1/(CC_1)$ we conclude that $\Pi(f)v\in W$ if $f\in O$ and
$p(v)<1$.

This proves the first statement. By taking inductive limits we infer
continuity of $\lim(\A_n(G)\times E)\to E^\omega$. For the
continuity of $\A(G)\times E\to E^\omega$ it now suffices to verify
that $\lim(\A_n(G)\times E)$ and $\A(G)\times E=(\lim \A_n(G))\times
E$ are isomorphic. The map
$$\lim(\A_n(G)\times E)\to (\lim \A_n(G))\times E$$
is clearly bijective and continuous. The left hand side is
{\it LF}, and the right hand side is a product of
ultrabornological spaces, hence itself ultrabornological.
It follows that the open mapping theorem can be applied
(see \cite{MV}, Theorem 24.30 and Remarks 24.15, 24.36).
\end{proof}

For later use we note that $\A(G)$ contains a Dirac sequence.

\begin{lem}\label{dirac}
The heat kernel $h_t$ belongs to $\A(G)$ for each $t>0$.
Let $E$ be an {\it F}-representation. Then
$\Pi(h_t)v\to v$ in $E$ for all $v\in E$.
\end{lem}

\begin{proof}
The convergence in $E$ is Nelson's theorem (see Section
\ref{sec::analyticreps}). The heat kernel belongs to $\A(G)$ for all
$t>0$ by \cite{GKL}, Thm.~4.2.
\end{proof}

\begin{rem}\label{diracrem}
It follows from the proof of \cite{GKL}, Thm.~4.2, that there exists
a common $m$ such that $h_t\in \A_m(G)$ for all $t>0$.
\end{rem}

\subsection{$\A(G)$-tempered representations}
As we have seen that there is continuous algebra action of $\A(G)$
on the analytic vectors of {\it F}-representations, we shall make
this property part of a definition.

\begin{dfn} A representation
$(\pi, E)$ is called $\A(G)$-{\rm tempered} if for all $f\in \A(G)$
and $v\in E$ the vector valued integral
$$\Pi(f)v=\int_G f(g) \pi(g) v\ dg$$
converges in $E$, and $(f,v)\mapsto\Pi(f)v$
defines a continuous algebra action $$\A(G)\times E\to E.$$
\end{dfn}

\begin{ex} \label{ex1}(a) For every F-representation $(\pi, E)$
both $(\pi, E)$ itself and $(\pi, E^\omega)$ are $\A(G)$-tempered
according to Proposition \ref{A action}. In particular this holds
for all Banach representations and also for $E=\cR(G)$  with the
left action (so that $E^\omega=\A(G)$).
\par \noindent (b) If $(\pi, E)$ is an  $\A(G)$-tempered
representation and $F\subset E$ is
a closed $G$-invariant subspace, then the induced representations on $F$ and $E/F$
are  $\A(G)$-tempered.
\end{ex}

\section{Analytic globalizations of Harish-Chandra
modules}\label{sec::globalize}

In this section we will assume that $G$ is a real reductive group.
Let us fix a maximal compact subgroup $K<G$. We say that a complex
vector space $V$ is a $(\gf, K)$-module if $V$ is endowed with a Lie
algebra action of $\gf$ and a locally finite group action of $K$
which are compatible in the sense that the derived and restricted
actions of $\kf$ agree and, in addition,
$$k \cdot (X\cdot v) = (\Ad(k)X) \cdot (k \cdot v) \qquad (k\in K, X\in\gf, v\in V)\, .$$
We call a $(\gf, K)$-module {\it admissible} if for any irreducible
representation $(\sigma, W)$ of $K$ the multiplicity space $\Hom_K
(W, V)$ is finite dimensional. Finally, an admissible $(\gf,
K)$-module is called a {\it Harish-Chandra module} if $V$ is
finitely generated as a $\U(\gf)$-module. Here, as usual, $\U(\gf)$
denotes the universal enveloping algebra of $\gf$.

\par By a {\it globalization} of a Harish-Chandra module $V$ we
understand
a representation $(\pi, E)$ of $G$ such that the space of $K$-finite vectors
$$E_K:=\{ v\in E\mid \dim \Span_\C\{\pi( K)v\} <\infty\}$$
is $(\gf, K)$-isomorphic to $V$ and dense in $E$. Density of $E_K$
is automatic whenever $E$ is quasi-complete, see \cite{HCdsII},
Lemma 4. Each element $v\in E$ allows an expansion in $K$-types
$v=\sum_{\tau\in \hat K} v_\tau$, where
$v_\tau=\dim\tau\,\pi(\chi_\tau)v\in E_K$. Here, the integral over
$K$ that defines $\pi(\chi_\tau)v$ may take place in the completion
of $E$, but $v_\tau$ belongs to $E_K$ by density and finite
dimensionality of $K$-type spaces.

A {\it Banach (F-, analytic, $\A(G)$-tempered) globalization} is a
globalization by a Banach ({\it F}-, analytic, $\A(G)$-tempered)
representation. Note that according to Harish-Chandra
\cite{HCrepBanach}, $E_K\subset E^\omega$ if $E$ is a Banach
globalization. In general, the orbit map of a vector $v\in E_K$ is
weakly analytic (see Remark \ref{weak anal}).

According to the subrepresentation theorem of
Casselman (see \cite{W} Thm.~3.8.3), $V$ admits a Banach-globalization $E$.
The space $E^\omega$ is then an analytic $\A(G)$-tempered globalization.

If $V$ is a Harish-Chandra module, we denote by $\tilde V$ the
Harish-Chandra module dual to $V$, i.e.~the space of $K$-finite
linear forms on $V$ (see \cite{W}, p.~115). We note that if $E$ is a
globalization of $V$, then $\tilde V$ embeds into $E'$ and
identifies with the subspace of $K$-finite continuous linear forms
(see \cite{C}, Prop.~2.2). Furthermore $\tilde V$  separates on $E$.
Since the matrix coefficients $x\mapsto \xi(\pi(x)v)$ for $v\in V,
\xi\in \tilde V$ are real analytic functions on $G$, they are
determined by their germs at $\1$. It follows that these functions
on $G$ are independent of the globalization (see \cite{C}, p.~396).

\subsection{Minimal analytic globalizations}

Let $V$ be a Harish-Chandra module and $\v=\{v_1, \ldots, v_k\}$ be
a set of $\U(\gf)$-generators.  We shall fix an arbitrary
$\A$-tempered globalization $(\pi, E)$ and regard $V$ as a subspace
in $E$.

On the product space $\A(G)^k=\A(G)\times\dots\times\A(G)$
with diagonal $G$-action, we consider the $G$-equivariant
map
$$\Phi_\v: \A(G)^k \to E, \ \ \f=(f_1, \ldots, f_k )\mapsto
\sum_{j=1}^k \Pi(f_j) v_j\ ,$$ and write $I_\v$ for its kernel. This
map is evidently continuous, and thus $I_\v$ is a closed
$G$-invariant subspace of $\A(G)^k$. We note that $\f\in I_\v$ if
and only if $\sum_j \int f_j(g) \xi(\pi(g)v_j)\,dg=0$ for all
$\xi\in \tilde V$. It follows that $I_\v$ is independent of the
choice of globalization. Furthermore, the dependence on generators
is easily described: If $\v'$ is another set of generators, say $k'$
in number, there exists a $k\times k'$-matrix $u$ of elements from
$\U(\gf)$ such that $\f\in I_\v$ if and only if $R_u\f\in I_{\v'}.$

Since $I_\v$ is closed and $G$-invariant, the quotient
$$V^\mg:= \A(G)^k / I_\v$$
carries a representation of $G$ which we denote by $(\pi, V^\mg)$.
It is independent of the choice of the globalization $(\pi,E)$ and
(up to equivalence) of the set $\v$ of generators.

\begin{lem}\label{minimality of globalization}
Let $V$ be a Harish-Chandra module.
Then the following assertions hold:
\begin{enumerate}
\item $V^\mg$ is an analytic $\A(G)$-tempered globalization of $V$.
\item $V^\mg=\Pi(\A(G))V$, that is, $V^\mg$ is spanned by the vectors of form $\Pi(f)v$.
\item If $(\lambda, F)$ is any $\A$-tempered globalization of $V$, then the identity mapping
$V\to F$ lifts to a $G$-equivariant continuous injection $V^\mg\to F^\omega$.
\end{enumerate}
\end{lem}

\begin{proof} (i) It follows from the definition that $V^\mg$
is analytic (see Lemma \ref{subrep}) and $\A(G)$-tempered (see
Example \ref{ex1}(b)). It remains to be seen that $(V^\mg)_K$ is
$(\gf,K)$-isomorphic to $V$. By definition, $\Phi_\v$ induces a
continuous $G$-equivariant injection $V^\mg\to E$. In particular
$(V^\mg)_K$ is isomorphic to a $(\gf,K)$-submodule of $V=E_K$.
Moreover as $\A(G)$ contains a Dirac sequence by Lemma \ref{dirac},
and as we may assume $E$ to be a Banach space, each generator $v_j$
belongs to the $E$-closure of the image of $V^\mg$. By admissibility
and finite dimensionality of $K$-types, $v_j$ belongs to $(V^\mg)_K$
for each $j$. Thus $(V^\mg)_K\simeq V$ and (i) follows. Assertions
(ii) and (iii) are clear.
\end{proof}

Because of property (iii), we shall refer to $V^\mg$ as the {\it
minimal $\A(G)$-tempered globalization} of $V$. We record the
following functorial properties of the construction.

\begin{lemma}\label{functorial}
Let $V,W$ be Harish-Chandra modules.
\begin{enumerate}
\item Every $(\gf,K)$-homomorphism
$T\colon W\to V$ lifts to a unique intertwining operator
$T^\mg\colon W^\mg\to V^\mg$ with restriction $T$ on $W=(W^\mg)_K$
and with closed image.
\item Assume that $W\subset V$ is a submodule. Then
\begin{enumerate}
\item $W^\mg$ is equivalent with a subrepresentation of $V^\mg$
on a closed invariant subspace,
\item $(V/W)^\mg$ is equivalent with the quotient representation
$V^\mg/W^\mg$.
\end{enumerate}
\end{enumerate}
\end{lemma}

\begin{proof} (i) Let $\tilde T\colon \tilde V\to\tilde W$ denote
the dual map of $T$, and observe that
$$\tilde T\xi(\pi(g)w)=\xi(\pi(g)Tw)$$
for all $w\in W,\xi\in \tilde V$ and $g\in G$. Indeed, these are
analytic functions of $g$ whose power series at $\1$ agree because
$T$ is a $\gf$-homomorphism. It follows that if we choose generators
$w_1,\dots,w_l$ for $W$ and $v_1,\dots,v_k$ for $V$ such that
$v_j=Tw_j$ for $j=1,\dots,l$, then the inclusion map
$\f\mapsto(\f,\0)$ of $\A(G)^l$ into $\A(G)^k$ takes $I_\w$ into
$I_\v$. Hence this inclusion map induces a map
$$T^\mg\colon \A(G)^l/I_\w\to \A(G)^k/I_\v$$
which is continuous, intertwining and has closed image.
Moreover, this map restricts to $T$ on $W$, since it maps
each generator $w_j$ to $v_j=Tw_j$.

(ii) is obtained from (i) with $T$ equal to
(a) the inclusion map $W\to V$ or (b) the quotient map $V\to V/W$.
\end{proof}

Our next concern will be to realize the analytic globalizations
inside of Banach modules.

\begin{prop} \label{p=Be} Let $(\pi, E)$ be an analytic
$\A(G)$-tempered globalization of a Harish-Chandra
module $V$. Then there exists a Banach representation $(\sigma, F)$ of $G$ and
a continuous $G$-equivariant injection $(\pi, E)\to (\sigma, F)$.
\end{prop}

\begin{proof}  We fix
generators $\xi=\{\xi_1, \ldots, \xi_l\}$ of the dual Harish-Chandra
module $\tilde V\subset E'$ and put
$U:=\{ v\in E\mid \max_{1\leq j \leq l} |\xi_j(v)|<1\}$. Then
$U$ is an open neighborhood of $0$ in $E$.

\par Fix $m\in \N$ such that $\A_m(G)$  contains a Dirac sequence
(see Remark \ref{diracrem}).
As $\A_m(G)\times E \to E$ is continuous,
we find an open neighborhood $O$ of $0$ in $\A_m(G)$
and an open neighborhood $W$ of $0$ in $E$ such that $\Pi(O)W\subset U$.
We may assume that $O$ is of the type $O=\{ f\in \A_m(G)\mid q(f)<1\}$
where
$$q(f)=\sup_{g\in G\atop z\in \Omega} |f(zg)| e^{N d(g)} $$
for some $N\in \N$ and $\Omega\subset V_m$ compact.
Define the normed space $X:=(\A_m(G), q)$. It follows from
Remark \ref{right action} that the right regular
action of $G$ is a representation by bounded operators on $X$.
Let $F:=(X^*)^l$ be the topological dual of $X^l $ and $\sigma$
the corresponding dual
diagonal action of $G$. Note that $F$ is a Banach space,
being the dual of a normed space, so that
$\sigma$ is a Banach representation. We claim that
the map
$$\phi: E\to F, \ \ v\mapsto \left(\f=(f_1, \ldots, f_l)\mapsto \sum_{j=1}^l \xi_j(\Pi(f_j)v)\right) $$
is $G$-equivariant, continuous and injective. Equivariance is clear,
and in order to establish continuity we fix a closed convex
neighborhood $\tilde O$ of $0$ in $F$. We may assume that $\tilde O$
is a polar  of the form $\tilde O=[B^l]^o$   where $B$ is a bounded
set $B\subset X$. Because $B$ is bounded, there exists $\lambda>0$
such that $B\subset \lambda O$. Choosing $\tilde W:=\frac1{\lambda}
W$ we have $\phi(\tilde W)\subset \tilde O$, as
\begin{align*} \phi(\tilde W)(B^l) &\subset \frac1{l}\phi(W)(O^l) \subset
\frac1{l} \sum_{j=1}^l \xi_j(\Pi(O)W)\\
& \subset \frac1{l} \sum_{j=1}^l \xi_j(U)\subset \{ z\in \C\mid
|z|\leq 1\}\, .\end{align*} It remains to be shown that $\phi$ is
injective. Suppose that $\phi(v)=0$. Then $\phi(v_\tau)=0$ for each
element $v_\tau$ in the $K$-finite expansion of $v$, so that we may
assume $v$ is $K$-finite. Then for all $f\in \A_m(G)$ and $\eta\in
\tilde V$ one would have $\eta(\Pi(f)v)=0$. Since $K$-finite matrix
coefficients are independent of globalizations, we conclude by Lemma
\ref{dirac} that $\eta(v)=0$ and hence $v=0$.
 \end{proof}

\subsection{The minimal analytic globalization of a spherical principal series representation}

Let $G=KAN$ be an Iwasawa decomposition of $G$ and denote by $M$
the  centralizer of $A$ in $K$, i.e. $M=Z_K(A)$.
Then $P=MAN$ is a minimal parabolic subgroup.
Let us denote by $\af,  \nf $ the Lie algebras of $A$ and $N$ and define $\rho\in \af^*$
by $\rho(X) =\frac 12\tr (\ad X|_\nf)$, $X\in \af$.
For $\lambda\in \af_\C^*$ and $a\in A$ we set $a^\lambda:=e^{\lambda(\log a)}$.

\par The {\it smooth spherical principal series
with parameter $\lambda\in \af_\C^*$} is defined by
$$V_\lambda^\infty:=\{ f\in C^\infty(G)\mid \forall man\in P\ \forall g\in G
:\ f(mang)=a^{\rho+\lambda} f(g)\}\, . $$ The action of $G$ on
$V_\lambda^\infty$ is by right displacements in the arguments, and
in this way we obtain a smooth {\it F}-representation $(\pi_\lambda,
V_\lambda^\infty)$ of $G$. We denote the Harish-Chandra module of
$V_\lambda^\infty$ by $V_\lambda$.

\par It is useful to observe that the restriction mapping to $K$,
$$\Res_K: V_\lambda^\infty \to C^\infty (M\backslash K),$$
is an $K$-equivariant isomorphism of Fr\'echet spaces, and
henceforth we will identify $V_\lambda^\infty$ with
$C^\infty(M\backslash K)$. The space $V_\lambda$  of $K$-finite
vectors in $V_\lambda^\infty$ is then identified as a $K$-module
with the space $C(M\backslash K)_K$ of $K$-finite functions on $M\bs
K$.

Likewise, the Hilbert space $\H_\lambda:=L^2(M\bs K)$ is provided
with the representation $\pi_\lambda$. The space of smooth vectors
for this representation is
$\H_\lambda^\infty=V^\infty_\lambda=C^\infty(M\bs K)$, and the space
of analytic vectors is the space
$\H_\lambda^\omega=V_\lambda^\omega:=C^\omega(M\bs K)$ of analytic
functions on $M\bs K$ with its usual topology.

\begin{thm} \label{th=ps}For every $\lambda\in \af_\C^*$ one has
$$\Pi_\lambda(\A(G))V_\lambda= C^\omega(M\bs K)\, .$$
In particular $V_\lambda^\mg\simeq V_\lambda^\omega=C^\omega(M\bs K)$
as analytic representations.
\end{thm}

The proof of this theorem is similar to the corresponding result in
the smooth case (see \cite{BK}, Section 4). Note that from Lemma
\ref{minimality of globalization} we have
$$\Pi_\lambda(\A(G))V_\lambda=V_\lambda^\mg \subset V_\lambda^\omega$$
with continuous inclusion. As the space $V_\lambda^\mg$ admits a web
(see \cite{MV}, 24.8 and 24.28) and $C^\omega(M\backslash K)$ is
ultrabornological (see \cite{MV}, 24.16), we can apply the open
mapping theorem (\cite{MV}, 24.30) to obtain an identity of
topological spaces from the set-theoretical identity. It thus
suffices to prove that for each $v\in V_\lambda^\omega$ there exists
$\xi\in V_\lambda$ and $F\in \A(G)$ such that $\Pi(F)\xi=v$.

We need some technical
preparations. Let us denote by $\gf=\kf +\pf$ the Cartan decomposition
of $\gf$, and write $\theta$ for the corresponding Cartan
involution. Let $(\cdot, \cdot)$ be a non--degenerate invariant
bilinear form on $\gf$ which is positive definite on $\pf$ and
negative definite on $\kf$. Then $\la\cdot, \cdot\ra = -(\theta
\cdot, \cdot)$ defines an inner product on $\gf$, which we use to
identify $\gf$ and $\gf^*$. We write $|\cdot|$ for the norms induced
on $\gf$ and $\gf^*$.

\par Let $X_1, \ldots, X_s$ be an orthonormal basis of $\kf$ and $Y_1, \ldots, Y_l $ be an orthonormal basis
of $\pf$. We define elements in the universal enveloping algebra $\U(\gf)$ by
$$\Delta=\sum_{j=1}^s X_j^2 +\sum_{i=1}^l Y_i^2\, , \qquad
\Delta_K=\sum_{j=1}^s X_j^2 \quad\hbox{and} \quad
\rC:=\Delta-2\Delta_K\, .$$ Note that $\rC$ is a Casimir element. In
particular, it belongs to the center of $\U(\gf)$.

\par Let $\tf\subset \kf$ be a maximal torus.
We fix a positive system of the root system $\Sigma(\tf_\C, \kf_\C)$
and identify the unitary dual $\hat K$ via their highest weights
with a subset of $i\tf^*$.  If $(\tau, W_\tau)$ is an irreducible
representation of $K$, then $\Delta_K$ acts as the scalar multiple
$|\tau +\rho_\kf|^2 -|\rho_\kf|^2$. For every $\tau\in \hat K$ we
denote by $\chi_\tau\in C(K)$ the normalized character
$\chi_\tau(k)=(\dim W_\tau)^{-1} \tr \tau(k)$. Note that $C(K)$ acts
on $\A(G)$ by left convolution.
\par We denote the left regular representation of $G$ on $\A(G)$ by $L$.
The following proposition will be crucial in the proof of
Theorem \ref{th=ps}.

\begin{prop}\label{p11} Let $(c_\tau)_{\tau\in \hat K}$ be a sequence of
complex numbers and $(a_\tau)_{\tau\in \hat K}$ a sequence of
elements in $G$. Assume that
$$ |c_\tau| \leq C e^{-\e|\tau|},\qquad
d(a_\tau)\leq c_1\log(1+|\tau|) + c_2  $$
for some $C,\e, c_1,c_2>0$.
Let $f\in \A(G)$. Then
$$F:=\sum_{\tau\in \hat K} c_\tau \chi_\tau * L(a_\tau)f  \in \A(G)\, .$$
\end{prop}
\begin{proof} As $(L, \cR(G))$ is an {\it F}-representation,
it follows from \cite{GKL} that
$h\in \cR(G)$ is belongs to $\A(G)$
if and only if  there exists an $M>0$ such that for all
$N\in \N$ there exists a constant $C_N>0$ with
\begin{equation} \label{ee} \sup_{g\in G} e^{Nd(g)}
|\Delta^k h(g)|\leq C_N M^{2k} (2k)!\, \end{equation}
for all $k\in \N$.
\par Observe that $\Delta=\rC + 2\Delta_K$. For every $h\in \cR(G)$ one has
\begin{equation} \label{e1} \Delta_K (\chi_\tau* h) =
(|\tau+\rho_\kf|^2 -|\rho_\kf|^2)  \chi_\tau* h\, .\end{equation}
Moreover as $\rC$ is central we obtain for every $g\in G$ and $h\in \A(G)$
that
\begin{equation} \label{e2}
\rC (\chi_\tau* L(g)h)= \chi_\tau * L(g)(\rC h)\, .\end{equation}
\par Let now $f\in \A(G)$. As $f$ is an analytic vector for $\cR(G)$,
hence also for $L^2(G)$, we find (see \cite{GW}, Cor. 4.4.6.4) a
constant $M_1>0$ such that for all $N>0$ there exists a constant
$C_N>0$ such that
\begin{equation} \label{est1}
\sup_{g\in G} e^{Nd(g)} |\rC^k f(g)|\leq C_N M_1^{2k} (2k)!\, . \end{equation}
\par We first estimate $\Delta^k (\chi_\tau*L(a_\tau)f)$. For that we employ
(\ref{e1}) and (\ref{e2}) in order to obtain that
\begin{align*} \Delta^k (\chi_\tau*L(a_\tau)f)&=\sum_{j=0}^k {k\choose j} \rC^j
(2\Delta_K)^{k-j}  (\chi_\tau*L(a_\tau)f)\\
&=\sum_{j=0}^k 2^{k-j} {k\choose j} (|\tau+\rho_\kf|^2 -|\rho_\kf|^2)^{k-j}
(\chi_\tau*L(a_\tau)\rC^jf)\, .
\end{align*}
For $N>0$ we thus obtain using (\ref{est1}) that
\begin{align*}\sup_{g\in G}\, & e^{N d(g)}  |\Delta^k (\chi_\tau*L(a_\tau)f)(g)|\\
&\leq C_N  2^{2k} \sum_{j=0}^k (1+ |\tau|)^{2(k-j)}
\cdot \sup_{g\in G} e^{N d(g)} |L(a_\tau)\rC^j f(g) | \\
&\leq C_N' 2^{2k}  e^{N d(a_\tau)} \sum_{j=0}^k M_1^{2j} (1+ |\tau|)^{2(k-j)}(2j)!\\
&\leq C_N''   M_2^{2k}  \sum_{j=0}^k   (1+ |\tau|)^{2(k-j)+Nc_1} (2j)!
\end{align*}
for some $C_N, M_2>0$ independent of $\tau$. Using these
inequalities for $F$ we arrive at
$$\sup_{g\in G}  e^{N d(g)} |\Delta^k F(g)|\leq
C_N''  M_2^{2k} \sum_{\tau\in \hat K} \sum_{j=0}^k |c_\tau|
(1+|\tau|)^{2(k-j)+c_1} (2j)!\, .$$ From the lemma below we obtain
that
$$\sum_{\tau\in \hat K} |c_\tau| (1+|\tau|)^{2(k-j)+c_1}
\leq C M^{2k-2j} (2k-2j)! $$
for some constants $C,M>0$ independent of $k,j$. Since
$$\sum_{j=0}^k  (2k-2j)!(2j)!\leq 2^{2k}(2k)!$$
we conclude that
$F$ satisfies the estimates (\ref{ee}).
\end{proof}

\begin{lemma}
Let $\epsilon>0$. There exist $C,M>0$ such that
$$\sum_{\tau\in \hat K} e^{-\epsilon|\tau|} (1+|\tau|)^n\leq C\,M^n\,n!$$
for all $n\in\N$.
\end{lemma}

\begin{proof} We assume for simplicity that $K$ is semisimple. The proof
is easily adapted to the general case. The set $\hat K$ is parametrized
by a semilattice in $i\tf^*$, say $\hat K=\{m_1\tau_1+\dots+m_l\tau_l\mid
m_1,\dots,m_l\in\N\}$. We shall perform the summation over $\hat K$ by
summing over $m\in\N$, and over those elements $\tau$ for which the maximal
$m_j$ is $m$. There are $l\,m^{l-1}$ such elements, and they all satisfy
$am\leq|\tau|\leq bm$ for some $a,b>0$ independent of $m$.
It follows that the sum above is dominated
by
$$\sum_{m\in\N} l m^{l-1} e^{-\epsilon am} (1+bm)^n.$$
The given estimate now follows easily.
\end{proof}

Before we give the proof of Theorem \ref{th=ps}, we
recall some harmonic analysis for the compact
homogeneous space $M\bs K$. We denote by
\def\hatKM{K_M^\wedge}$\hatKM$
the $M$-spherical part of $\hat K$, that is,
the equivalence classes of irreducible representations $\tau$
for which the space $V_\tau^M$ of $M$-fixed vectors
is non--zero. Then
\begin{equation}\label{Peter-Weyl}
L^2(M\bs K)=
\hat\oplus_{\tau\in\hatKM} \Hom(V_\tau^M,V_\tau)
\end{equation} by the Peter--Weyl
theorem. We write $v=\sum_\tau v_\tau$ for the corresponding
decomposition of a function $v$ on $M\bs K$ and note that with the
right action of $k\in K$ on $L^2(M\bs K)$ we have
$[\pi(k)v]_\tau=\tau(k)\circ v_\tau$.

Furthermore,
$$C^\omega(M\bs K)=\big\{v=\sum_\tau v_\tau
\big|\, \exists \epsilon, C>0\, \forall\tau: ||v_\tau||\leq
Ce^{-\epsilon|\tau|} \big\}\, ,$$ where $\|v_\tau\|$ denotes the
operator norm of $v_\tau$.

Let $\tau\in\hatKM$. The integral
$\delta_\tau(k)=\dim(\tau)\int_M \chi_\tau(mk)\,dm$
of the character is bi-invariant
under $M$. The components of $\delta_\tau$
in the decomposition (\ref{Peter-Weyl}) are all $0$ except
the $\tau$-component,
which is the inclusion operator $I_\tau$ of
$V_\tau^M$ into $V_\tau$.

\begin{proof}
We can now finally give the proof of Theorem \ref{th=ps}.
Let $v=\sum_\tau v_\tau\in C^\omega(M\bs K)$
be given, and let $\epsilon>0$ be as above.

It follows from \cite{BK}, Section 6, that there exists a
$K$--finite function $\xi\in V_\lambda$, and for each
$\tau\in\hatKM$ elements $a_\tau\in A$ and $c_\tau\in\C$ such that
$$d(a_\tau)\leq c_1\log(1+|\tau|)+c_2,\quad
|c_\tau|\leq 2(1+|\tau|)^{c_3}
$$
for some constants $c_1,c_2,c_3>0$ independent of $\tau$,
and such that
$$
R_\tau:=\delta_\tau-c_\tau[\pi_\lambda(a_\tau)\xi]_\tau
$$
satisfies
$\|R_\tau\|\leq 1/2
$
for all $\tau$. By integration of $\xi$ and $R_\tau$ over $M$
we can arrange that they are both $M$--biinvariant.

We now choose a function $f\in\A(G)$ such that $\Pi(f)\xi=\xi$. It
exists because $\Pi(\sum_{\tau\in F}\chi_\tau*h_t)\xi$ converges to
$\xi$ for $t\to 0$ and some finite set $F$ of $K$-types by Lemma
\ref{dirac}, so that $\xi$ belongs to the closure of a finite
dimensional subspace of $\Pi(\A(G))\xi$. According to Proposition
\ref{p11}, the function
$$F=\sum_{\tau}
c_\tau e^{-\frac12\epsilon|\tau|}\chi_\tau*L(a_\tau)f$$
belongs to $\A(G)$. An easy calculation shows that
$$\Pi(F)\xi=\sum_\tau e^{-\frac12\epsilon|\tau|}(\delta_\tau-R_\tau)\ .$$

Being of type $\tau$ and $M$--biinvariant, $R_\tau$ corresponds
in (\ref{Peter-Weyl}) to an operator
$R_\tau\in \End(V_\tau^M)$.
Since $\|R_\tau\|\leq 1/2$, the operator
$I_\tau-R_\tau\in \End(V_\tau^M)$ is invertible
with $\|(I_\tau-R_\tau)^{-1}\|\leq 2$.
Then $v_\tau(I_\tau-R_\tau)^{-1}\in\Hom(V_\tau^M,V_\tau)$
with $\|v_\tau(I_\tau-R_\tau)^{-1}\|\leq 2Ce^{-\epsilon|\tau|}$.
It follows that the function on $M\bs K$
with the expansion
$$\sum_{\tau} e^{\frac12\epsilon|\tau|} v_\tau(I_\tau-R_\tau)^{-1}$$
belongs to $C^\omega(M\bs K)$. We denote by $h(k^{-1})$
this function, so that $h$ is a right $M$--invariant function on $K$.
Another easy calculation now shows that
$$\Pi(h)\Pi(F)\xi=
\sum_\tau v_\tau=v\, ,$$ and hence $h*F\in \A(G)$ is the function we
seek.
\end{proof}

\subsection{Unique analytic globalization}
The goal of this section is to prove the following
version of Schmid's minimal globalization theorem
(\cite{KaS}, Theorem 2.13).

\begin{thm}\label{uniqueglobalization}
Let $V$ be a Harish-Chandra module. Every analytic $\A(G)$-tempered
globalization of $V$ is isomorphic to $V^\mg$.

In particular, if $(\pi,E)$ is an arbitrary
{\it F}-globalization of $V$, then $$E^\omega\simeq V^\mg.$$
\end{thm}

\begin{proof} We first treat the case of an irreducible
Harish-Chandra module~$V$.

We first claim that $V$ admits a Hilbert globalization $\H$ such
that $\H^\omega=\Pi(\A(G))V$, and hence in particular (see Lemma
\ref{minimality of globalization} (ii))
$$\H^\omega\simeq V^\mg.$$
\par In case $V=V_\lambda$ we can take $\H_\lambda=L^2(M\backslash K)$ and the assertion follows
from Theorem \ref{th=ps}. If the Harish-Chandra module is of the type $V=V_\lambda\otimes W$ where
$W$ is a finite dimensional $G$-module, then $\H=\H_\lambda\otimes W$ is a Hilbert globalization with
$\H^\omega=\H_\lambda^\omega\otimes W$. A straightforward generalization of \cite{BK}, Lemma 5.4, yields that
$$(\Pi_\lambda\otimes \Sigma)(\A(G))V=\H^\omega\, .$$
Finally, every irreducible Harish-Chandra module is a subquotient of
some $V_\lambda\otimes W$ (see for example \cite{LW}, Thm.\ 4.10),
and the claim follows by Lemma \ref{functorial}.

\par
Let now $(\pi, E)$ be an arbitrary analytic $\A(G)$-tempered
globalization of $V$. We aim to prove $E\simeq V^\mg$. From Lemma
\ref{minimality of globalization} we know that $V^\mg$ injects
$G$-equivariantly and continuously
into $E=E^\omega$, hence it suffices to establish
surjectivity of the injection.

\par We now fix the Hilbert globalization $\H$ of above. In view of
Proposition \ref{p=Be} we can embed $(\pi, E)$ into a Banach
globalization $F$ of $V$. As $E$ is analytic, we obtain a continuous
$G$-equivariant injection $E\to F^\omega$.  In order to proceed we
recall the Casselman-Wallach theorem (cf. \cite{C}, \cite{W}, or
\cite{BK} for a more recent proof) which implies that $F^\infty$ is
equivalent to $\H^\infty$ as {\it F}-representation. It follows, see
Corollary \ref{analytic-smooth}, that $F^\omega\simeq \H^\omega$.
Collecting the established isomorphisms, we have
$$V^\mg\to E\subset F^\omega\simeq \H^\omega\simeq V^\mg.$$
The surjectivity follows from the completeness of $\H^\omega$ (see
Proposition \ref{p=c}).

\par Finally, we prove the case of an arbitrary Harish-Chandra module.
As Harish-Chandra modules
have finite composition series, it suffices to prove the following
statement: Let $0\to V_1\to V\to V_2\to 0$ be an exact sequence of
Harish-Chandra modules and suppose that both $V_1$ and $V_2$ have
unique analytic $\A(G)$-tempered globalizations. Then so does $V$.

Let $E$ be an analytic $\A(G)$-tempered globalization of $V$. Let
$E_1$ be the closure of $V_1$ in $E$ and $E_2=E/E_1$. Then $E_1$ and
$E_2$ are analytic $\A(G)$-tempered globalizations of $V_1$ and
$V_2$. By assumption we get $E_1=V_1^\mg$ and $E_2=V_2^\mg$, and
from Lemma \ref{functorial} we infer $V_2^\mg=V^\mg/V_1^\mg$.
Observe that in an exact sequence of topological vector spaces $0\to
E_1\to E\to E_2\to 0$ the topology on $E$ is uniquely determined by
the topology of $E_1$ and $E_2$ (see \cite{DS}, Lemma 1). We thus
conclude that $E=V^\mg$.
\end{proof}

We conclude by summarizing the topological properties of $V^\mg$.
Recall that an inductive limit $E = \lim_{n \to \infty} E_n$ of
Fr\'{e}chet spaces is called regular if every bounded set is
contained and bounded in one of the steps $E_n$.

\begin{cor} The minimal globalization $V^\mg$ is a
nuclear, regular, reflexive and complete inductive limit of
Fr\'{e}chet--Montel spaces.
\end{cor}
\begin{proof} Theorem \ref{uniqueglobalization} and Proposition \ref{p=c} imply that $V^\mg$ is complete.
Furthermore, it then follows from \cite{wen} and \cite{kuc} that
$V^\mg$ is regular and reflexive (see also Appendix B). It is an
inductive limit of Fr\'{e}chet--Montel spaces, because $\A(G)$ is an
inductive limit of Fr\'{e}chet--Schwarz spaces, 
and Hausdorff quotients of such spaces are Fr\'{e}chet--Montel. 
Nuclearity is inherited from $C^\omega(M\bs K)$, which is the 
strong dual of a nuclear Fr\'{e}chet space,
and this property is preserved when
passing to the quotient of a finite dimensional tensor product.
Finally, a Fr\'echet space is nuclear if and only
if its strong dual is nuclear (see \cite{J}, Section 21.5).
\end{proof}

\section{Appendix A. Vector-valued holomorphy}\label{AppHol}

Here we collect some results about analytic functions with
values in a locally convex Hausdorff topological vector space $E$.
Let $\Omega\subset\C^n$ be open.

It is a natural and common assumption that $E$ is sequentially complete.
Let us recall that under this assumption an $E$-valued function
$f$ on $\Omega$ is said to be holomorphic if it
satisfies one of the following conditions, which are
equivalent in this case:

\par\smallskip
(a) $f$ is weakly holomorphic, that is, the scalar function
$z\mapsto \zeta(f(z))$ is holomorphic for each continuous linear form
$\zeta\in E'$;

(b)  $f$ is $\C$-differentiable in each variable at each $z\in\Omega$;

(c)  $f$ is infinitely often $\C$-differentiable at each $z\in\Omega$;

(d)  $f$ is continuous and is represented by a converging power series
expansion with coefficients in $E$, in a neighborhood of each $z\in\Omega$.
\par\smallskip

In general, the conditions (c) and (d) are mutually equivalent
and they imply (a) and (b). This follows by regarding
$f$ as a function into the completion $\bar E$ of $E$
(see \cite{Glockner}, Prop. 2.4).
We shall call a function $f: \Omega\to E$
{\it holomorphic} if (c) or (d) is satisfied,
or equivalently, if it is holomorphic into $\bar E$ with
$E$-valued derivatives up to all orders.

Let $M$ be an $n$-dimensional complex manifold. An $E$-valued
function on $M$ is called holomorphic if all its coordinate
expressions are holomorphic. We denote by $\O(M, E)$ the space of
$E$-valued holomorphic functions on $M$. Endowed with the compact
open topology, it is a Hausdorff topological vector space, which is
complete whenever $E$ is complete.

The following isomorphism of topological
vector spaces is useful.

\begin{lemma}\label{lemma: productiso}
Let $M$ and $N$
be complex manifolds, then
\begin{equation}\label{productiso}
\O(M\times N, E)\simeq \O(M,\O(N,E))
\end{equation}
under the natural map $f\mapsto (x\mapsto f(x,\,\cdot\,))$
from left to right.
\end{lemma}

\begin{proof} Apart from the statement that
$x\mapsto f(x,\,\cdot\,)\in\O(N,E)$ is holomorphic, this is
straightforward from definitions. It is clear that
$f(x,\,\cdot\,)\in\O(N,E)$. By regarding $\O(N,E)$ as a subspace of
$\O(N,\bar E)$ and noting that it carries the relative topology, we
reduce to the case that $E$ is complete, so that condition (b)
applies. Assume for simplicity that $M=\C$. What needs to be
established is then only that the complex differentiation
$$\frac{\partial f}{\partial x}(x,y)=\lim_{h\to 0}
\frac1h[f(x+h,y)-f(x,y)]\in E$$
is valid locally uniformly with respect to $y\in N$.
This follows from uniform continuity on compacta of the derivative.
\end{proof}

\section{Appendix B. Topological Properties of $\A(G)$}\label{AppTop}

While the topology of a general inductive limit of Fr\'{e}chet
spaces may be complicated, $\A(G)$ inherits certain properties from
the steps $\A(G)_n$.

\begin{thm} \label{thm::Acomplete}
The algebra $\A(G)$ is regular, complete and reflexive.
\end{thm}

A regular inductive limit of Fr\'{e}chet--Montel spaces is known to
be reflexive \cite{kuc} and complete \cite{wen}, so that we only
have to show regularity. The following criterion from \cite{wen},
Theorem 3.3, in terms of interpolation inequalities will be
convenient:

\begin{prop}
An inductive limit $E = \lim_{n \to \infty} E_n$ of
Fr\'{e}chet-Montel spaces is regular if and only if for some
fundamental system $\{p_{n,\nu}\}_{\nu\in\N}$ of seminorms on $E_n$:
$\forall n \ \exists m>n\ \exists \nu\ \forall k>m\ \forall \mu\
\exists \kappa\ \exists C\ \forall f \in E_n$
\begin{equation}\label{eq::wroth}
p_{m,\mu}(f)\ \leq\ C(p_{k,\kappa}(f) +p_{n,\nu}(f))\ .
\end{equation}
\end{prop}

In the case of $\A(G)$, condition (\ref{eq::wroth}) should be
thought of as a weighted geometric relative of Hadamard's Three
Lines Theorem. To verify it, we need to introduce some notions from
complex and Riemannian geometry, starting with the appropriate
differential operators.

By common practice we identify the Lie algebra $\mathfrak{g}_\C$ with the
space of right--invariant vector fields on $G_\C$, where $X \in
\mathfrak{g}_\C$ corresponds to the differential operator
$$\widetilde{X} u(x) = \frac{d}{dt}\Big |_{t=0} u(\exp(-tX)x)
\quad (x \in G_\C, \ u \in C^\infty(G_\C)). $$

If we denote the complex structure on the Lie algebra
$\mathfrak{g}_\C$ by $J$, the Cauchy--Riemann operators
$\overline{\partial}_Z$ and $\partial_Z$ associated to $Z \in
\mathfrak{g}_\C$ are given by $\overline{\partial}_{Z} :=
\widetilde{Z} + i \widetilde{JZ}$ resp.~$\partial_{Z} :=
\widetilde{Z} - i \widetilde{JZ}$.

In this section it will be convenient to replace the
left $G$-invariant metric $\bfg$ on $G$ used in Section
\ref{subsec::growth}
by a right invariant one,
which we shall denote by the same symbol. Note that the corresponding
distance functions $d$ on $G$ are equivalent (see (\ref{exponentialestimateweight}).
The function
$$K(\exp(JX)g) := \frac{1}{2}|X|^2:=\frac{1}{2}\mathbf{g}_\1(X,X)$$
endows a sufficiently small complex neighborhood $VG$ of $G$ with a
right $G$--invariant K\"{a}hler structure. To see this, choose an
orthonormal basis $\{X_j\}_{j=1}^l$ of $\mathfrak{g}$ with respect
to the metric. A straightforward computation results in
$$\partial_{X_i}\overline{\partial}_{X_j} K(\1) = \mathbf{g}_\1(X_i, X_j),$$
so that the complex Hessian $(Z_1, Z_2)\mapsto
\partial_{Z_1}\overline{\partial}_{Z_2} K(\1)$ defines a positive
definite Hermitian form on $\mathfrak{g}_\C$. By continuity and
invariance, positivity extends to give a K\"{a}hler metric on a
small neighborhood $VG$.

The complex Laplacian $$\Delta_\C = \sum_{j=1}^l
\partial_{X_j}\overline{\partial}_{X_j}=\sum_{j=1}^l\widetilde{X_j}^2+\widetilde{JX_j}^2,$$
agrees with the K\"{a}hler Laplacian up to first--order terms and
maps real--valued functions to real--valued functions. Therefore the
following weak maximum principle holds:

\begin{lem}\label{lem::maxprinc}
If $u\in C^2(VG)$ is real--valued with a local maximum in $z \in
VG$, then $$\Delta_\C u(z)\leq 0.$$
\end{lem}

As $\Delta_\C$ is a trace of the complex Hessian, we may rely on
well--known results about plurisubharmonic functions to conclude:

\begin{lem}\label{lem::subharm}
For $u\in \mathcal{O}(VG)$, $\Delta_\C u = 0$ and $\Delta_\C \log|u|
\geq 0$.
\end{lem}

So while it may be less obvious how to control applications of
$\Delta_\C$ to the Riemannian distance function $d$ on $G$,
$\Delta_\C$ annihilates the holomorphically regularized distance
function $\tilde{d}:=e^{-\Delta_\mathbf{g}} d$ from \cite{GKL}. This
is going to be useful in the proof of Theorem \ref{thm::Acomplete},
and the following Lemma, which is  shown as in Lemma 4.3 of \cite{GKL}
collects the key properties of $\tilde{d}$.

\begin{lem}\label{lem::dtilde}
a) The function $\tilde{d}$ extends to a function in
$\mathcal{O}(UG)$ for some neighborhood $U$ of $\1 \in
G_\C$. \\
b) For all $U'\Subset U$, $\sup_{zg \in U'G}
|\tilde{d}(zg)-d(g)|<\infty$ and $\widetilde{X_j}\ \tilde{d}$ as
well as $\widetilde{JX_j}\ \tilde{d}$ are bounded on $U'G$ for all
$j$.
\end{lem}

Before finally coming to the proof of Theorem \ref{thm::Acomplete},
we introduce an equivalent representation of $\A(G)$ based on
geometrically more convenient neighborhoods.
If we define for $n \in \N$, $\nu
\in \N_0$, the neighborhoods
$$\widetilde{V}_n :=\left\{\exp(JX)\in G_\C \mid
|X| <\frac{1}{n}\right\},$$
$$\Omega_{n}^{\nu}:=\left\{\exp(JX)\in G_\C \mid
|X|<\frac{1}{n+(\nu+2)^{-1}}\right\}$$ and associated subspaces
of $\A(G)$,
\begin{align*}
&\widetilde{\A(G)}_n := \\&\qquad\Big\{f \in \mathcal{O}(\widetilde{V}_n
G) \mid \forall \nu \in \N : p_{n,\nu}(f):= \sup_{g\in G, z\in
\Omega_{n}^{\nu}} |f(zg)|\ e^{\nu d(g)}<\infty \Big\},
\end{align*}
then $\A(G)$ is again an inductive limit $\lim_{n \to
\infty}\widetilde{\A(G)}_n$ of Fr\'{e}chet--Montel spaces. Condition
(\ref{eq::wroth}) translates into
\begin{equation}\label{eq::wrothA}
\sup_{zg\in \Omega_{m}^{\mu}G} |f(zg)|\ e^{\mu d(g)} \leq C\
(\sup_{zg\in \Omega_{k}^{\kappa}G} |f(zg)|\ e^{\kappa d(g)} +
\sup_{zg \in \Omega_{n}^{\nu}G} |f(zg)|\ e^{\nu d(g)})\
\end{equation}
for $f \in \widetilde{\A(G)}_n$.

To show this, let $n$ sufficiently large, $0 \not \equiv f \in
\widetilde{\A(G)}_n$, $m=n+1$, $\nu= 0$, $k>m$ and $\mu\in \N$,
and consider
$$u(z) = \log|f(z)|
+ N(z)\ D(z)$$ on $\widetilde{V}_nG\setminus
\widetilde{V}_{k+1}G$, where we choose $N(\exp(JX)g) =
N(\exp(JX)) = \bar{\nu}(
|X|^{-2\alpha}-(n+\frac{1}{2})^{2\alpha})$ and $D(z) =
D_0+\mathrm{Re}\ \tilde{d}(z)$ for some $\bar{\nu},\alpha, D_0>0$.
First note that $\Delta_\C u >0$ if $D_0$ and $\alpha$ are
sufficiently large. Indeed, by Lemma \ref{lem::subharm} it is enough
to show $\Delta_\C (N(z) D(z))>0$. But $\Delta_\C D=0$, so that
\begin{align*} \Delta_\C (N(z) D(z)) & = (\Delta_\C N(z))\ D(z) +\\
& +2 \sum_{j=1}^l \left\{\widetilde{X_j}N(z) \widetilde{X_j} D(z) +
\widetilde{JX_j}N(z) \widetilde{JX_j} D(z)\right\}\, .
\end{align*}
With $D \geq 1$ on $\widetilde{V}_nG$ for large $D_0$ by Lemma
\ref{lem::dtilde}, we only have to show that $$\Delta_\C N(z)>
\overline{D}\ \max_{j=1,\dots,l} \{|\widetilde{X_j}N(z)|,
|\widetilde{JX_j}N(z)|\}$$ on $\widetilde{V}_nG$ for large $n$ and
$\overline{D}=2\ \sup \{|\widetilde{X_j}D|,\ |\widetilde{JX_j}D| :
j=1,\dots, l\}$. By $G$--invariance, it is sufficient to do so in
$z=\exp(\varepsilon JX)$ close to $\varepsilon=0$. The
Baker--Campbell--Hausdorff formula implies
\begin{align*} \exp(tJX_j)\exp(\varepsilon JX) &= \exp(\varepsilon JX+tJX_j +
\mathcal{O}(\varepsilon t^2)+\mathcal{O}(\varepsilon^2
t))\cdot \\
&\cdot \exp(\frac{1}{2}\varepsilon t [JX_j,JX])\ ,
\end{align*}
 so that
\begin{eqnarray*}
\widetilde{JX_j}\ N(\exp(\varepsilon JX)) &=&
\frac{d}{dt}\Big|_{t=0}N(\varepsilon JX+tJX_j
+\mathcal{O}(\varepsilon t^2)+
\mathcal{O}(\varepsilon^2 t))\\
&=&-2\alpha\bar{\nu}\
\varepsilon^{-1-2\alpha}\ \frac{ \mathbf{g}_\1(X_j,X)}{\mathbf{g}_\1(X,X)^{\alpha+1}}+
\mathcal{O}(\varepsilon^{-2\alpha}),\\
\end{eqnarray*}
Similarly,
\begin{align*}
&(\widetilde{JX_j})^2 \ N(\exp(\varepsilon JX))\\&\qquad=
2\alpha\bar{\nu}\ \varepsilon^{-2-2\alpha}\
\frac{2(\alpha+1)\mathbf{g}_\1( X_j,X)^2-\mathbf{g}_\1(X_j,X_j)
\mathbf{g}_\1( X,X)}{\mathbf{g}_\1(X,X)^{\alpha+2}}
\end{align*}
up to terms of
order $\varepsilon^{-1-2\alpha}$. Summing over $j$ establishes the
assertion for large $\alpha$ and small $\varepsilon$, hence for
large $n$.

For $\kappa \geq 0$, set $S_{n} := \sup_{\partial \Omega_{n}^{0} G}
u$ and $S_{k}^\kappa := \sup_{\partial \Omega_{k}^{\kappa} G} u$.
Because $u(z)$ is bounded from above and $\leq \max\{S_{k}^\kappa,
S_{n}\}$ on $\partial \Omega_{k}^{\kappa}G \cup
\partial \Omega_{n}^{0}G$, the maximum principle, Lemma
\ref{lem::maxprinc}, assures $$u(z) \leq \max\{S_{k}^\kappa,
S_{n}\}$$ in $\Omega_{n}^{0}G\setminus \Omega_{k}^\kappa G$, or
\begin{eqnarray*}
|f(z)|\ e^{N(z) D(z)} &\leq &e^{\max\{S_{k}^\kappa, S_{n}\}}
\\&\leq&\sup_{w\in \partial \Omega_{k}^\kappa G} |f(w)|\ e^{N(w) D(w)}+ \sup_{w\in \partial \Omega_{n}^{0}G} |f(w)|\ e^{N(w)
D(w)}.
\end{eqnarray*}

As $\widetilde{V}_m \Subset \Omega_{n}^{0}$, we may choose
$\bar{\nu}$ such that $N|_{\widetilde{V}_mG\setminus
\widetilde{V}_{k+1}G} \geq \mu$. Setting $\kappa :=
\sup_{\widetilde{V}_kG\setminus \widetilde{V}_{k+1}G} N \geq \mu$
we obtain
\begin{eqnarray*}
\sup_{z\in \Omega_{m}^\mu G} |f(z)|\ e^{\mu D(z)} &\leq & \sup_{z\in
\Omega_k^\kappa G} |f(z)|\ e^{\kappa D(z)}+\sup_{z\in
\partial
\Omega_n^0 G} |f(z)|\\
&\leq & \sup_{z\in \Omega_k^\kappa G} |f(z)|\ e^{\kappa
D(z)}+\sup_{z\in \Omega_n^0 G} |f(z)|.
\end{eqnarray*}
Lemma \ref{lem::dtilde} implies $d(z)-C \leq D(z) \leq d(z)+C$ for
some $C>0$, and Theorem \ref{thm::Acomplete} follows.

\begin{rem}
It would be interesting to better understand the topology of
$\A(G)^N/I$ for a stepwise closed, $\A(G)$--invariant subspace $I$.
Because $\widetilde{\A(G)}_{n}$ is even Fr\'{e}chet--Schwarz, the
quotients $\widetilde{\A(G)}_{n}^N/(I\cap\widetilde{\A(G)}_{n}^N)$
are Fr\'{e}chet--Montel and one might hope to verify condition
(\ref{eq::wroth}) as before. However, adapting the above proof
requires strong assumptions on $I$, and general Hausdorff quotients
$\A(G)^N/I$ are likely to be incomplete: For a convex domain $\Omega
\subset \R^n$, the space of test functions $\mathcal{D}(\Omega)$ is
isomorphic to a similar weighted space of holomorphic functions by
Paley--Wiener's theorem. However, given any non--surjective
differential operator $A$ on $\mathcal{D}'(\Omega)$, the quotient of
$\mathcal{D}(\Omega)$ by the image of $A^t$ will be incomplete.
\end{rem}

\end{document}


\title[Analytic representations]
{Corrigendum: ``Analytic representation theory  of Lie groups: General theory and
analytic globalizations of Harish--Chandra modules''\\
$[$
Compos. Math. 147 (2011), 1581--1607 $]$}

\author[Gimperlein]{Heiko Gimperlein}
\author[Kr\"otz]{Bernhard Kr\"otz}
\author[Schlichtkrull]{Henrik Schlichtkrull}

\address{Maxwell Institute for Mathematical Sciences and Department of Mathematics, Heriot--Watt University \\ Edinburgh, EH14 4AS \\ \newline and Institut f\"ur Mathematik, Universit\"at Paderborn\\Warburger Str. 100 
\\ 33098 Paderborn \\ email: h.gimperlein@hw.ac.uk}
\address{Institut f\"ur Mathematik, Universit\"at Paderborn\\Warburger Str. 100  
\\ 33098 Paderborn \\email: bkroetz@gmx.de}
\address{Department of Mathematics\\
University of Copenhagen\\Universitetsparken 5 \\
DK-2100 Copenhagen
\\email: schlicht@math.ku.dk}

\date{\today}
\thanks{}
\begin{abstract}
Some errors in the original article are corrected.\end{abstract}

\maketitle
For a representation $\pi$ of a connected Lie group $G$ on a topological vector space $E$
we defined in \cite{GKS} a vector subspace $E^\omega$ of $E$ of analytic vectors. Further we equipped 
$E^\omega$ with an inductive limit topology.  We called a representation $(\pi, E)$ {\it 
analytic}  if $E=E^\omega$ as topological vector spaces.

\par Some mistakes in the paper have been pointed out 
by Helge Gl\"ockner (see \cite{G}).
For a representation $(\pi, E)$ and a closed $G$-invariant subspace $F$ of $E$ we asserted 
in Lemma 3.6 (i) that $F^\omega=E^\omega\cap F$ as a topological space.
Based on that we further asserted in Lemma 3.6 (ii) that the inclusion 
$E^\omega/F^\omega \to (E/F)^\omega$ is continuous and in Lemma 3.11
that if $(\pi,E)$ is analytic then so is the restriction to $F$. 
However, there is a gap in the proof of the first assertion,
and presently it is not clear to us whether 
the above statements are then true in this generality (for unitary representations $(\pi, E)$ they are straightforward).
Our proof does give the following weaker version of the two lemmas:

\bigskip\noindent
{\bf Lemma 1.} {\it Let $(\pi,E)$ be a representation, and let $F\subset E$ be a closed
invariant subspace. Then
\begin{enumerate}
\item[(i)] $F^\omega=E^\omega\cap F$ as vector spaces and with 
continuous inclusion $F^\omega\to E^\omega$.
\item[(ii)] $E^\omega/E^\omega\cap F \subset (E/F)^\omega$ continuously.
\item[(iii)] If $(\pi,E)$ is an analytic representation, then $\pi$ induces an analytic
representation on  $E/F$.
\end{enumerate}}
Indeed, for (iii) note that if $E$ is analytic, $E/F = E^\omega/E^\omega\cap F \subset (E/F)^\omega$ continuously by (ii), and $(E/F)^\omega \subset E/F$ continuously.

\bigskip Further we asserted in Proposition 3.7 
a general completeness property of the functor 
which associates $E^\omega$ to $E$. However, there is a gap in the proof, 
which asserts that $v_i\to v$ in the topology of $E^\omega$. As statements in this generality 
are not needed for the main result, we can leave out the proposition (together with Remark 3.8).

\par Attached to $G$ we introduced a certain analytic convolution algebra $\A(G)$. A central theme 
of the paper is the relation of analytic  representations of $G$ to algebra representations
of $\A(G)$ on $E$: $\A(G)\times E \to E$.
In Proposition 4.2 (ii) we claimed that the bilinear map $\A(G)\times \A(G) \to \A(G)$ is 
continuous. However, the proof shows only separate continuity. For a similar reason we need to weaken Proposition 4.6 to:

\bigskip\noindent {\bf Proposition 2.} {\it 
Let $(\pi, E)$ be an F-representation.  
The assignment  
$$(f, v)\mapsto \Pi(f)v:=\int_G f(g)\pi(g)v\,dg$$
defines a continuous bilinear map 
$$\A_n(G)\times E\to E_n $$
for every $n\in \N$, and a separately continuous map 
$$\A(G)\times E\to E^\omega$$
(with convergence of the defining integral in $E^\omega$).
Moreover, if $(\pi, E)$ is a Banach representation, then
the latter bilinear map is continuous.}

\begin{proof} The first statement is proved in the article, and thus
only the statement for $\pi$ a Banach representation remains to 
be proved. We repeat the 
first part of the proof, 
now with $p$ denoting the fixed norm of $E$. The constants
$c,C$ such that
$$p(\pi(g)v)\leq C e^{c d(g)} p(v) \qquad (g\in G, v\in E)$$
and $N,C_1$ such that
$$C_1:=\int_G e^{(c-N)d(g)}\,dg<\infty,$$
are then all fixed, and so is $\epsilon=1/(CC_1)$.

Let $n\in\N$ and an open $0$-neighborhood
$W_n\subset E_n$ be given. 
We may assume
$$W_n=\{ v\in E_n\mid p(\pi(K_n)v)<\epsilon_n\}$$
with $K_n\subset GV_n$ compact and $\epsilon_n>0$. 
Let
$$O_n:=\{ f\in\O(V_nG)\mid \sup_{z\in K_n, g\in G} |f(z^{-1}g)| e^{Nd(g)} <
\epsilon\epsilon_n\}
\subset \A_n(G).$$
The computation in the given proof shows that if
$f\in O_n$ and $p(v)<1$
then $\Pi(f) v\in W_n$.
The asserted bi-continuity of $\A(G) \times E \to E^\omega$ follows. 
\end{proof}

As a consequence we obtain as in Example 4.10(a), but only for
Banach representations $(\pi,E)$, that
$E^\omega$ is $\A(G)$-tempered. 
In particular $\A(G)$ need not itself be
$\A(G)$-tempered, and we need to 
replace Lemma 5.1(i) by the following weaker version: 

\bigskip\noindent {\bf Lemma 3.} {\it $V^\mg$ is an analytic globalization of $V$ and it 
carries an algebra action 
$$(f,v)\mapsto \Pi(f)v,\qquad \A(G)\times V^\mg\to V^\mg,$$
of $\A(G)$, which is separately continuous.}

\bigskip

The main result of the paper Theorem 5.7, has two
statements concerning a Harish-Chandra module $V$ with a
globalization $E$:

\begin{enumerate}{\it
 \item If $E$ is analytic $\A(G)$-tempered then  $E=V^\mg$
 \item If $E$ is an $F$-globalization then $E^\omega=V^\mg$. }
\end{enumerate}

The proof, which relied on Lemma 3.11 and Proposition 4.6
respectively, needs to be corrected. The proof of (1) if
$V$ is irreducible needs no modification. For the general
case it can be adjusted as follows.

Like in the paper, it suffices to consider an exact sequence of Harish--Chandra modules $0\to V_1 \to V \to V_2 \to 0$, where both $V_1$ and $V_2$ have unique analytic $\A(G)$-tempered globalizations. We show that the same holds for $V$.

Let $E_1$ be the closure of $V_1$ in $E$ and $E_2 = E/E_1$. By Lemma 1(iii), $E_2$ is an analytic $\mathcal{A}(G)$-tempered globalization of $V_2$, so that by assumption $E_2 = V_2^{\mg} = \mathcal{A}(G) V_2$ as topological vector spaces.

In a first step we prove that $E_1 = V_1^\mg = \mathcal{A}(G) V_1$ as vector spaces. For that we note first that $E_1$ is ${\mathcal A}(G)$-tempered
and that $V_1^\mg\subset E_1$ continuously. 
Next, by Proposition 5.3 (which holds for any ${\mathcal A}(G)$-tempered representation), we may embed $E_1 \subset F_1$ continuously into a Banach globalization of $F_1$ of $V_1$. Moreover, the proof shows that the embedding is compatible with the action by $\mathcal{A}(G)$. It follows that 
$E_1^\omega \subset F_1^\omega$ continuously and as $\A(G)$-modules.
 Further note that since $E$ is analytic, from Lemma 1(i) we also obtain $E_1^\omega = E^\omega \cap E_1 = E_1$ as vector spaces. Hence $V_1^\mg \subset E_1 \subset F_1^\omega$. 
By assumption $V_1$ has unique $\A(G)$-tempered globalization, hence $F_1^\omega \simeq V_1^\mg$. 
Therefore $V_1^\mg \subset E_1 \subset F_1^\omega \simeq V_1^\mg$. 
As these maps respect the structure as $\mathcal{A}(G)$--module, the inclusion is also surjective: 
$V_1^\mg=E_1$.

Being an inductive limit, $E_1 = F_1^\omega$ is an ultrabornological space, and 
$V_1^\mg$ is webbed (see the reference in the proof of Proposition 4.6). We conclude from
the open mapping theorem that $V_1^\mg = E_1$ also as topological vector spaces. 

With Lemma 5.2 we now have a diagram of topological vector spaces
$$
\begin{matrix}
0&\to&V_1^\mg&\to &V^\mg&\to&V_2^\mg&\to&0\\
&&\parallel&&\downarrow&&\parallel&\\
0&\to&E_1&\to&E&\to&E_2&\to&0\\
 \end{matrix}
$$
where the vertical arrow in the middle signifies the continuous
inclusion $V^\mg= \mathcal{A}(G) V  \subset E$, and where the rows are exact.
The five--lemma implies $V^\mg= E$ as a vector space, and as in the article we 
conclude from [DS79] that this is then a topological identity.

\par Finally, for (2) we recall from Corollary 3.5 that $(E^\infty)^\omega =E^\omega$. The Casselman-Wallach smooth globalization 
theorem asserts the existence of a Banach globalization $F$ of $V$ such that $F^\infty=E^\infty$ and 
therefore $F^\omega=E^\omega$. In particular,   $E^\omega$ is $\A(G)$-tempered  by Proposition 2. 
Now (1) applies. 

\bigskip 
{\bf Acknowledgement.}
The authors wish to thank Helge Gl\"ockner for pointing out the discussed 
mistakes.